\documentclass[12pt]{amsart}
\usepackage{amssymb,latexsym}
\usepackage{enumerate}
\usepackage{amscd}
\usepackage{verbatim}

\theoremstyle{plain}
\newtheorem{theorem}{Theorem}
\newtheorem{corollary}[theorem]{Corollary}
\newtheorem{lemma}[theorem]{Lemma}
\newtheorem{proposition}[theorem]{Proposition}

\theoremstyle{definition}

\theoremstyle{remark}

\numberwithin{equation}{section}
\numberwithin{theorem}{section}
\numberwithin{definition}{section}

\def\pair#1#2{\langle #1, #2\rangle}

\def\Mof#1{\Cal M(A)}

\def\Set{{\text{\bf Set}}}

\def\Unif{{\text{\bf Unif}}}
\def\HUnif{{\text{\bf HUnif}}}
\def\CHUnif{{\text{\bf CHUnif}}}
\def\CU{{\text{\bf CU}}}
\def\HCU{{\text{\bf HCU}}}
\def\CHCU{{\text{\bf CHCU}}}

\newcommand{\Con}{\operatorname{Con}}

\newcommand{\Div}{\operatorname{Div}}

\newcommand{\Fg}{\operatorname{Fg}}

\newcommand{\OpSemiUnif}{\operatorname{SemiUnif}}

\newcommand{\MyIm}{\operatorname{Im}}
\newcommand{\nat}{\operatorname{nat}}

\newcommand{\Fil}{\operatorname{Fil}}

\newcommand{\I}{\operatorname{I}}
\newcommand{\OpUnif}{\operatorname{Unif}}

\def\congruence{on\-gru\-ence\discretionary{-}{}{-}}
\def\conM/{c\congruence mod\-u\-lar}
\def\ConM/{C\congruence mod\-u\-lar}
\def\conD/{c\congruence dis\-trib\-u\-tive}
\def\ConD/{C\congruence dis\-trib\-u\-tive}
\def\conP/{c\congruence per\-mut\-a\-ble}
\def\ConP/{C\congruence per\-mut\-a\-ble}
\def\conMity/{\conM/\-i\-ty}
\def\ConMity/{\ConM/\-i\-ty}
\def\conDity/{c\congruence dis\-trib\-u\-tiv\-i\-ty}
\def\ConDity/{C\congruence dis\-trib\-u\-tiv\-i\-ty}
\def\conPity/{c\congruence per\-mut\-a\-bil\-i\-ty}
\def\ConPity/{C\congruence per\-mut\-a\-bil\-i\-ty}
\def\usprv/{un\-der\-ly\-ing-set-pre\-ser\-ving}

\def\ie/{{i.e.}}
\def\Ie/{{I.e.}}
\def\eg/{{e.g.}}
\def\Eg/{{E.g.}}
\def\etc/{{etc.}}

\newdimen\mysubdimen
\newbox\mysubbox

\def\subwhat#1#2#3{{
\setbox\mysubbox=\hbox{#3}
\mysubdimen=\wd\mysubbox
\setbox\mysubbox=\hbox{$#1#2$}
\ifnum\mysubdimen>\wd\mysubbox
\vtop{
\hbox to\mysubdimen{\hfil\box\mysubbox\hfil}
\nointerlineskip
\hbox{#3}}
\else
\mysubdimen=\wd\mysubbox
\vtop{
\box\mysubbox
\nointerlineskip
\hbox to\mysubdimen{#3}}
\fi
}}

\def\mydate{\ifcase\month Month0 \or January \or
February \or March \or April \or May \or June \or July \or
August \or September \or October \or November \or
December \fi \number\day, \number\year}

\begin{document}

\title[Algebras with a Compatible Uniformity]{Algebras
with a Compatible Uniformity}
\author{William H. Rowan}
\address{PO Box 20161 \\
         Oakland, California 94620}
\email{whrowan@member.ams.org}


\date{\today}

\begin{abstract}
Given a variety of algebras $\mathbf V$, we study
categories of algebras in $\mathbf V$ with a compatible
structure of uniform space.  The lattice of compatible
uniformities of an algebra, $\OpUnif A$, can be
considered a generalization of the lattice of congruences
$\Con A$. Mal'cev properties of $\mathbf V$ influence the
structure of $\OpUnif A$, much as they do that
of $\Con A$. The
category $\mathbf V[\CHUnif]$ of complete, Hausdorff
such algebras in the variety $\mathbf V$ is
particularly interesting; it has a
factorization system $\pair{\mathbf E}{\mathbf M}$, and
$\mathbf V$ embeds into $\mathbf V[\CHUnif]$ in such a
way that $\mathbf E\cap\mathbf V$ is the subcategory of
onto and $\mathbf M\cap\mathbf V$ the subcategory of
one-one homomorphisms.
\end{abstract}
\maketitle

\section*{Introduction}

In this paper, we study algebras
having a compatible structure of uniform space. We will
call these \emph{uniform algebras} in this paper. We study
the lattice of compatible uniformities
$\OpUnif A$ on an algebra $A$, and show that Mal'cev
properties on a variety of algebras $\mathbf V$ influence
the structure of $\OpUnif A$ for $A\in\mathbf V$. We also
define a number of categories of uniform algebras in
$\mathbf V$.  We are particularly interested in the
category $\mathbf V[\CHUnif]$ of complete, Hausdorff
uniform algebras in $\mathbf V$.

We warn the reader that the term \emph{uniform algebra},
although used in this paper, has been used previously, in
a different sense, for referring to algebras such that the
lattice of subalgebras, or the lattice of congruences,
has the property that any two elements which are not
$\bot$ have a meet which is also not $\bot$.

We see $\OpUnif A$ as a
generalization of the congruence lattice $\Con A$. $\Con
A$ embeds into
$\OpUnif A$ by a lattice homomorphism that preserves
arbitrary joins. Also, given $\mathcal U\in\OpUnif A$, we
can complete $A$ with respect to $\mathcal U$, yielding a
uniform algebra we denote by $A/\mathcal U$ because this
construction generalizes the formation of a quotient
algebra by a congruence of $A$.

Mal'cev properties on a variety $\mathbf V$ influence the
structure of $\OpUnif A$ for $A\in\mathbf V$, much as
they do that of $\Con A$. We show that if $\mathbf V$ is
Mal'cev, then compatible uniformities on $A$ permute, in
a sense we define, and that $\OpUnif A$ is modular. We
show furthermore that if $\mathbf V$ is an arithmetical
variety, then $\OpUnif A$ is distributive. Finally, we
show that if $\mathbf V$ is congruence-modular, and thus
has a sequence of Day terms, then $\OpUnif A$ satisfies a
weakened form of the modular law.

The category $\mathbf V[\CHUnif]$ of complete, Hausdorff
uniform algebras in $\mathbf V$ can be considered as a
generalization of the category of algebras in $\mathbf
V$ (which we also denote by $\mathbf V$), in the following
sense: $\mathbf V$ embeds into $\mathbf V[\CHUnif]$, and
$\mathbf V[\CHUnif]$ admits a factorization system
$\pair{\mathbf E}{\mathbf M}$ such that $\mathbf
E\cap\mathbf V$ is the subcategory of onto and $\mathbf
M\cap\mathbf V$ the subcategory of one-one
homomorphisms.

The paper first defines uniform algebras in sections
\ref{S:Unif} to \ref{S:CongUnif}. Sections 
\ref{S:Joins} to \ref{S:CM} study the lattice of
compatible uniformities $\OpUnif A$, and the influence of
Mal'cev properties on its structure. Sections
\ref{S:Hausdorff} to \ref{S:Completion} study the
important operations of Hausdorffization and completion.
Section \ref{S:Limits} discusses limits and colimits in a
number of categories of uniform algebras. Finally,
sections \ref{S:FactSys} and \ref{S:Factorization} review
the theory of factorization systems in a category, and
discuss the important properties of the category $\mathbf
V[\CHUnif]$ that we have mentioned.

In this paper, the term \emph{proposition} has been
used for those results which are both elementary
consequences of the basic definitions, and well-known. We
have used \emph{theorem} for other results, even though
parts of some of them are straightforward and fairly
standard.

\section*{Preliminaries}

\subsection{Category theory} We follow
\cite{Mac L} in terminology and notation.

\subsection{Universal algebra} We assume the basic
definitions of Universal Algebra, as found, for example,
in
\cite{B-S}. However, in the definition of an algebra, we
prefer to allow an algebra to be empty.

\subsection{Lattices and Filters}
We use the notations $\bot$ and $\top$ for the least and
greatest elements of a lattice, assuming they exist. We
denote the interval sublattice from $a$ to $b$, in the
lattice $L$, by $\I_L[a,b]$.

A \emph{filter} in a lattice $L$ is a nonempty subset
$F\subseteq L$ which is closed under finite meets, and
closed upward, i.e., $x\leq y$ and $x\in F$ implies $y\in
F$. If we have a subset $B$ of a filter $F$, such that
$F$ consists of all elements of $L$ greater than or equal
to some element of $B$, then we say that $B$ is a
\emph{base} for $L$. A set $B\subseteq L$ is a base for a
filter of $L$ iff for any two elements $x$, $y\in B$,
there is an element $z\in B$ such that $z\leq x\wedge y$.

If $S$ is a subset of a lattice $L$, we
denote by $\Fg^L S$, or simply $\Fg S$, the \emph{filter
generated by $S$}. This is the smallest filter of $L$
containing $S$, and is the filter having the set of finite
meets of elements of $S$ as a base.

Filters in $L$ form a complete lattice $\Fil L$, with join
the intersection of filters and with meet of a tuple of
filters $\{\,F_i\,\}_{i\in I}$ equal to
$\Fg(\bigcup_iF_i)$.

\section{Binary Relations} \label{S:Relations}
In this paper, we
will make extensive use of binary relations on a set $S$.
The set of such relations forms a lattice, with
meet the intersection and join the union of
relations.

\subsection{Relational products}
Given two relations $U$, $V$ on
$S$, we will denote by $U\circ V$ the relational product
of $U$ and $V$, i.e., the relation $\{\,\pair
uv\mid\exists w\in S\text{ such that }u\mathrel
Uw\text{ and }w\mathrel Vv\,\}$. This is an associative
operation on relations on $S$.
If $\mathcal U_1$ and $\mathcal U_2$ are sets of
relations on $S$, then we denote by 
$\mathcal U_1\circ\mathcal U_2$ the set $\{\,U_1\circ
U_2\mid U_1\in\mathcal U_1, U_2\in\mathcal U_2\,\}$.
This is an associative operation on sets of relations on
$S$.

Given a relation $U$ on $S$ (or, a set $\mathcal U$ of
relations on $S$),
we will denote by $U^n$ (by $\mathcal U^n$) the $n$-fold
relational power, defined inductively by $U^1=U$,
$U^{k+1}=U\circ U^k$ (respectively, by $\mathcal
U^1=\mathcal U$, $\mathcal U^{k+1}=\mathcal U\circ
\mathcal U^k$).
We also define
$U^{-1}=\{\,\pair ts\mid\pair st\in U\,\}$, and $\mathcal
U^{-1}=\{\,U^{-1}\mid U\in\mathcal U\,\}$.

\subsection{Relations, functions, and operations}
If $f:S\to T$ is a function, and $U$ is a
binary relation on $S$, then we denote by $f(U)$ the set
of pairs
$\pair{f(x)}{f(y)}$ such that $\pair xy\in U$.
More generally,
if $U$ is a relation on $S$, and $\omega$ an
$n$-ary operation on
$S$, then we denote by $\omega(U)$ the set of pairs
$\pair{\omega(\vec a)}{\omega(\vec b)}$ such that
$a_i\mathrel Ub_i$ for all $i$. If $\mathcal U$ is a set
of binary relations on $S$, and $\omega$ is an operation,
then we denote by $\omega(\mathcal U)$ the set of
$\omega(U)$ for $U\in\mathcal U$.

\begin{lemma} \label{T:Composition} Let $S$ be a set,
$U$ and $V$ binary relations on $S$, and $\omega$ an
$n$-ary operation on $S$. Then we have
\begin{enumerate}
\item[(1)] $\omega(U\cap
V)\subseteq\omega(U)\cap\omega(V)$, and
\item[(2)]
$\omega(U\circ V)\subseteq\omega(U)\circ\omega(V).$
\end{enumerate}
\end{lemma}

If $f:S\to T$ is a function, and $V$ is a binary relation
on $T$, then we denote by $f^{-1}(V)$ the set of pairs
$\pair xy$ of elements of $S$ such that
$\pair{f(x)}{f(y)}\in V$. If $\mathcal V$ is a set of
binary relations, then we denote by $f^{-1}(\mathcal V)$
the set of $f^{-1}(V)$ for $V\in\mathcal V$.

If $\omega$ is an $n$-ary operation on $S$
for $n>0$, and $U$ is a binary relation on $S$, then we
denote by
$U^\omega$ the set of pairs
$\pair a{a'}$ such that $\omega(a,\vec b\,)\mathrel
U\omega(a',\vec b\,)$ for all $(n-1)$-tuples $\vec b$.

\begin{lemma}\label{T:Subset} If $\{\,U_i\,\}_{i\in I}$ is
a tuple of binary relations on $A$,
and
$\omega$ is an
$n$-ary operation for $n>0$, then
$\bigcup_iU_i^\omega\subseteq(\bigcup_iU_i)^\omega$.
\end{lemma}

\begin{proof}
We have
\begin{align*}
a\mathrel{\bigcup_iU_i^\omega} a'
&\implies a\mathrel{U_i^\omega}a' & &\text{for some $i$}
\\
&\implies \omega(a,\vec b\,)\mathrel{U_i}\omega(a',\vec
b\,) & &\text{for all $\vec b$, for some $i$} \\
&\implies \omega(a,\vec
b\,)\mathrel{(\bigcup_iU_i)}\omega(a',\vec b\,) &
&\text{for all $\vec b$}\\
&\implies
a\mathrel{(\bigcup_iU_i)^\omega}a'
\end{align*}
\end{proof}

\subsection{Compatible filters of relations}
Let $A$ be an algebra. Recall that if $U$ is a relation
on $A$, then we say that $A$ is \emph{compatible} if for
all basic operations $\omega$ on $A$, $\omega(U)\subseteq
U$. We will now generalize  this concept to filters
of relations.

\begin{lemma}
Let $A$ be an algebra, and $\mathcal U$ a filter of
relations on $A$. The following are equivalent:
\begin{enumerate}
\item[(1)] For each basic operation
$\omega$, and each
$U\in\mathcal U$, there exists $\bar U\in\mathcal U$ such
that $\omega(\bar U)\subseteq U$.

\item[(2)] For each term operation $t$, and each
$U\in\mathcal U$, there exists $\bar U\in\mathcal U$ such
that $t(\bar U)\subseteq U$.

\item[(3)] For each basic operation $\omega$,
$\Fg(\omega(\mathcal U))\leq\mathcal U$.

\item[(4)] For each term operation $t$, $\Fg(t(\mathcal
U))\leq\mathcal U$.

\end{enumerate}
\end{lemma}

\begin{proof} First, note that by
lemma~\ref{T:Composition}, $\omega(\mathcal U)$ and
$t(\mathcal U)$ are bases for filters of relations. Thus,
it is clear that (1)
$\iff$ (3) and (2)
$\iff$ (4). Clearly also, (2) $\implies$ (1); to show (1)
$\implies$ (2), it suffices to show that the set $\mathbf
T$ of operations
$t$ (of all arities) such that for all $U\in\mathcal U$,
$\exists\bar U\in\mathcal U$ with $t(\bar U)\subseteq U$,
is a clone. Clearly, the
$i^{\text{th}}$ of $n$ projection $\pi_{in}\in\mathbf T$
for all
$n$ and $i$; given
$n'$-ary $\omega'\in\mathbf T$ and an $n'$-tuple
$\vec\omega$ of $n$-ary elements of $\mathbf T$ and given
$U\in\mathcal U$, $\exists\bar U\in\mathcal U$ with
$\omega'(\bar U)\subseteq U$, and $\exists\hat
U\in\mathcal U$ such that $\omega_i(\hat U)\in\bar U$ for
all $i$, using the fact that $\mathcal U$ is a filter.
Then $(\omega'\vec\omega)(\hat U)\subseteq U$. Thus,
$\mathbf T$ is a clone.
\end{proof}

 Given a
filter
$\mathcal U$ of relations on
$A$, we say that $\mathcal U$ is a \emph{compatible
filter of relations} on $A$ if the equivalent conditions
of the lemma are satisfied.

\begin{theorem}
\label{T:CompClosure}
Let $A$ be an algebra. Then we have
\begin{enumerate}
\item[(1)] If $\{\,\mathcal U_i\,\}_{i\in I}$ is a tuple
of compatible filters of relations on $A$, then
$\bigwedge_i\mathcal U_i$ is compatible;
\item[(2)] if $\mathcal U$ and $\mathcal V$ are compatible
filters of relations on $A$, then so is $\Fg(\mathcal
U\circ\mathcal V)$; and
\item[(3)] if $\mathcal U$ is a compatible filter of
relations on $A$, then so is $\mathcal U^{-1}$.
\end{enumerate}
\end{theorem}

\section{Uniform Spaces}\label{S:Unif}

In this section, we review the notion of a uniform space,
following \cite{Ito} in many respects; see also
the references cited there, including \cite{Is}.

\subsection{Uniformities and semiuniformities}
Let $S$ be a set, provided with a set
$\mathcal U$ of binary relations on $S$ satisfying
\begin{enumerate}
\item[(U1)] if $U\in\mathcal U$, and $U\subseteq V$, then
$V\in\mathcal U$;
\item[(U2)] if $U$, $V\in\mathcal U$, then $U\cap
V\in\mathcal U$;
\item[(U3)] if $U\in\mathcal U$, then $\Delta_S\subseteq
U$, where $\Delta_S$ is the diagonal $\{\,\pair ss:s\in
S\,\}$;
\item[(U4)] if $U\in\mathcal U$, then $U^{-1}\in\mathcal
U$;
\item[(U5)] if $U\in\mathcal U$, then $V\circ V\subseteq
U$ for some $V\in\mathcal U$.
\end{enumerate}
Then $\mathcal U$ is called a \emph{uniformity} on $S$,
and
$\pair S{\mathcal U}$, or simply $S$, is called a
\emph{uniform space}.

Of course, (U1) and (U2) simply state that $\mathcal U$ is
a filter of binary relations, and (U3) states that all of
the relations in $\mathcal U$ are reflexive.

If $\mathcal U$ is a set of relations on $S$, satisfying
(U1) through (U4) but not necessarily (U5), then we call
$\mathcal U$ a \emph{semiuniformity}. Semiuniformities
are not as important as uniformities but will be
useful technically.

If $\mathcal U_1$ and $\mathcal U_2$ are uniformities
(or semiuniformities) on
$S$, such that $\mathcal U_1\subseteq\mathcal U_2$, then
we say that
$\mathcal U_1$ is \emph{weaker} than $\mathcal U_2$, and
that
$\mathcal U_2$ is \emph{stronger} than $\mathcal U_1$.
That is, strength is the same as
inclusion. However, we prefer to order uniformities the
same way we order filters, by reverse inclusion, so we
will write
$\mathcal U_1\leq\mathcal U_2$ if
$\mathcal U_1$ is stronger than $\mathcal U_2$.

The strongest uniformity on $S$, the
\emph{discrete uniformity} containing all reflexive
relations on $S$, is the least element of $\OpUnif S$ and
$\OpSemiUnif S$, and the weakest uniformity, the
\emph{indiscrete uniformity} containing only $S\times S$,
is the greatest element.

\subsection{The lattices $\OpUnif S$ and $\OpSemiUnif S$}
We denote the set of uniformities on
$S$ by
$\OpUnif S$, and the set of semiuniformities by
$\OpSemiUnif S$.

\begin{theorem}
$\OpUnif S$ is a complete meet
subsemilattice of the lattice of filters of reflexive
relations on $S$.
\end{theorem}

\begin{proof}
Given a tuple $\mathcal U_i$ of uniformities,
the meet of the tuple in the lattice of filters of
reflexive relations is the filter $\mathcal U$ having as base
the set of finite intersections of elements of
$\bigcup_i\mathcal U_i$. $\mathcal U$ clearly satisfies (U1)
through (U4). To prove (U5), it suffices to show that if
$U$ is a finite meet of elements of $\bigcup_i\mathcal U_i$,
then there is a $V\in\mathcal U$ such that $V\circ V\subseteq
U$. Thus, let $U=U_1\cap\ldots\cap U_n$, where
$U_j\in\mathcal U_{i_j}$. By (U5) for the $\mathcal
U_{i_j}$, there are
$V_j\in\mathcal U_{i_j}$ such that $V_j\circ V_j\subseteq
U_j$. Then $V=V_1\cap\ldots\cap V_n\in\mathcal U$ and
$V\circ V\subseteq U$.
\end{proof}

It follows from the theorem that $\OpUnif S$, being a
complete meet semilattice, is also a complete lattice. We
will be examining various sets of uniformities, which form
complete meet subsemilattices of the lattice of
uniformities of the set in question. The corresponding
join operations can be problematic and will be studied in
sections~\ref{S:Joins} and
\ref{S:Malcev}.

The join operation in $\OpSemiUnif S$ is much more
tractable:

\begin{theorem} $\OpSemiUnif S$ is a complete
sublattice of the lattice of filters of binary relations
on $S$.
\end{theorem}

\begin{corollary} $\OpSemiUnif S$ is distributive.
\end{corollary}

\subsection{Bases for uniformities}
If $\mathcal U\in\OpUnif S$, then a base for $\mathcal U$,
 viewed simply as a filter in the set of binary
relations on $S$, is called a \emph{base for the
uniformity}
$\mathcal U$.

For example, it is easy to see that the symmetric elements
of $\mathcal U$ (those $U\in\mathcal U$ for which
$U^{-1}=U$) form a base for $\mathcal U$.

It is well known that a set $\mathcal B$ of binary
relations on a set $S$ is a base for a uniformity on $S$
iff the following conditions hold:
\begin{enumerate}
\item[(B2)] If $U$, $V\in\mathcal B$, then there is a
$W\in\mathcal B$ such that $W\subseteq U\cap V$;
\item[(B3)] if $U\in\mathcal B$, then $\Delta_S\subseteq
U$;
\item[(B4)] if $U\in\mathcal B$, then there is a
$V\in\mathcal B$ such that $V\subseteq U^{-1}$; and
\item[(B5)] if $U\in\mathcal B$, then there is a
$V\in\mathcal B$ such that $V\circ V\subseteq U$.
\end{enumerate}

\subsection{Uniform continuity}
If $\pair S{\mathcal U}$ and $\pair T{\mathcal V}$ are uniform
spaces, and $f:S\to T$ is a function, we say that $f$ is
\emph{ uniformly continuous} if for every $V\in\mathcal
V$, there is a
$U\in\mathcal U$ such that $f(U)\subseteq V$.

If $f:S\to T$ is a function from a set $S$ into a uniform
space $\pair T{\mathcal V}$, then the set of relations
$\{\,f^{-1}(V):V\in\mathcal V\,\}$ is a base
for a uniformity $f^{-1}(\mathcal V)=$ on $S$, called
the \emph{inverse image} of $\mathcal V$ under $f$.

On the other hand, if $f:S\to T$ and $\mathcal U\in\OpUnif S$,
then the meet of all uniformities $\mathcal V$ on $T$, such
that $f$ is then uniformly continuous from $\pair
S{\mathcal U}$ to $\pair T{\mathcal V}$, is a uniformity
$f_*(\mathcal U)$, which we call the
\emph{direct image} of $\mathcal U$ under $f$.

The following well-known proposition relates these three
concepts:

 \begin{proposition} Let $\pair S{\mathcal U}$ and
$\pair T{\mathcal V}$ be uniform spaces, and $f:S\to
T$ a function. Then the following are
equivalent:
\begin{enumerate}
\item[(1)] $f$ is uniformly continuous.
\item[(2)] $\mathcal U$ is stronger than the inverse image
uniformity $f^{-1}(\mathcal V)$.
\item[(3)] $\mathcal V$ is weaker than the direct image
uniformity $f_*(\mathcal U)$.
\end{enumerate}
\end{proposition}

In some important cases, we have a good description for
the direct image uniformity:

\begin{theorem}\label{T:DirectImage} Let $f:S\to T$ be a
function with kernel equivalence relation $\psi$, and let
$\mathcal U\in\OpUnif S$. Let $\mathcal V$ be the direct
image uniformity
$f_*(\mathcal U)$. If for all $U\in\mathcal U$, there is a $\bar
U\in\mathcal U$ such that $\bar U\circ\psi\circ\bar U\subseteq
\psi\circ U\circ\psi$, then the set
$\mathcal B=\{\,V_U\subseteq T^2:U\in\mathcal U\,\}$ where
$V_U=\Delta_T\cup f(U)$, is a base for $\mathcal V$. If
in addition $f$ is onto, then
$\mathcal V=\mathcal B$.
\end{theorem}

\begin{proof}
First, we must verify (B2) through (B5).

(B2): Let $U$, $U'\in\mathcal U$. Define $U''=U\cap U'$.
Then $V_{U''}\subseteq V_U\cap V_{U'}$.

(B3) is trivial.

(B4): For all $U$, $(V_U)^{-1}=V_{U^{-1}}$.

(B5): If $U\in\mathcal U$, let $\bar U\in\mathcal U$
be such that
$\bar U\circ\psi\circ\bar U\subseteq\psi\circ U\circ\psi$.
Then
\[
V_{\bar U}\circ V_{\bar U}=
V_{\bar U\circ\psi\circ\bar U}\subseteq
V_{\psi\circ U\circ\psi}=V_U.
\]

Thus, $\mathcal B$ is a base for a uniformity, and it is
clear that that uniformity is $\mathcal V$.

 If $f$ is
onto, then
$V_U=f(U)$ for all
$U$ and if
$V_U\subseteq V'$, $U\subseteq f^{-1}(V')$. Thus
$f^{-1}(V')\in\mathcal U$ and we have $V'=V_{f^{-1}(V')}$.
\end{proof}

\subsection{The category $\Unif$}
It is clear that the composition of two composable
uniformly continuous functions is uniformly continuous,
and that the identity function on a uniform space is
uniformly continuous. Thus, the class of uniform spaces,
and the uniformly continuous functions between them, form
a category, which we denote by $\Unif$.

\subsection{Uniformities and topology}
Every uniform space $\pair S{\mathcal U}$ has a natural
topology, where the closure of a subset $S'\subseteq S$
is the set of all $s\in S$ such that for all $U\in\mathcal U$, $\pair s{s'}\in U$ for some $s'\in S'$.

\begin{theorem}\label{T:Dense}
Let $\pair S{\mathcal U}$ be a uniform
space, and let $\mathcal V\in\OpUnif S$ be such that
$\mathcal U\leq\mathcal V$. Let $T\subseteq S$ be a
dense
subset of $S$ with respect to $\mathcal U$. Then
\begin{enumerate}
\item[(1)] $\mathcal V$ has a base of sets of the form
\[
U\circ V|_T\circ U
\]
where $U\in\mathcal U$ and $V\in\mathcal V$, and
$V|_T=\{\,\pair t{t'}\mid t,\,t'\in T\text{ and }t\mathrel
Vt'\,\}$.

\item[(2)] $\mathcal V$ is determined by $\mathcal V|_T$
as the filter generated by sets of the above form, where
$\mathcal V|_T=\{\,V_T\mid V\in\mathcal V\,\}$.
\end{enumerate}
\end{theorem}

\begin{proof}
Given $V\in\mathcal V$, let $\tilde V\in\mathcal V$
be symmetric and such that $\tilde V^5\subseteq V$. If
$s\mathrel{\tilde V}s'$, then because $T$ is dense in
$\pair S{\mathcal U}$ and $\tilde V\in\mathcal U$,
there
exist 
$t$, $t'\in T$ such that $t\mathrel{\tilde V}s$ and
$t'\mathrel{\tilde V}s'$. We have
$s\mathrel{\tilde V}t\mathrel{\tilde V}s\mathrel{\tilde V}
s'\mathrel{\tilde V}t'\mathrel{\tilde V}s'$. Thus,
\[
\tilde V\subseteq\tilde V\circ(\tilde V^3)|_T\circ\tilde
V\subseteq V;
\]
(1) follows.

To prove (2), we must show that every set of the form
$U\circ V|_T\circ U$ is an element of $\mathcal V$. Given
such an element, let $\tilde V\in\mathcal V$ be
symmetric and such that $\tilde V^3\subseteq V$ and let
$\bar U\subseteq U$ be a symmetric element of
$\mathcal U$ such that $\bar U\subseteq\tilde V$. Then as
before, we have
\[
\tilde V\subseteq\bar U\circ(\tilde V\circ\tilde
V\circ\tilde V)|_T\circ\bar U\subseteq U\circ V|_T\circ
U;
\]
(2) follows.
\end{proof}

\begin{corollary} \label{T:Dense2}
Let $\pair S{\mathcal U}$ be a uniform space, and $f:T\to
S$ a function from another set $T$ such that the image in
$S$ is dense. If $\mathcal V$, $\mathcal V'\in\OpUnif S$
are such that $\mathcal U\leq\mathcal V$ and $\mathcal
U\leq\mathcal V'$, and $f^{-1}(\mathcal
V)=f^{-1}(\mathcal V')$, then $\mathcal V=\mathcal V'$.
\end{corollary}

\begin{proof} Since $f^{-1}(\mathcal V)=f^{-1}(\mathcal
V')$, $\mathcal V|_{\MyIm f}=\mathcal V'|_{\MyIm f}$. By
part (2) of the theorem, this implies $\mathcal V=\mathcal
V'$.
\end{proof}

\section{Uniform Algebras}\label{S:UnifAlg}

\subsection{Finite powers of a uniform space}
In order to discuss uniformly continuous
$n$-ary operations on a uniform space
$\pair S{\mathcal
U}$, we must first describe the $n^{\text{th}}$
power of
$\pair S{\mathcal U}$ for $n$ finite. (We will discuss
arbitrary powers,
and more generally, limits and colimits, in
section~\ref{S:Limits}.)

\begin{theorem}
The $n^{\text{th}}$ power of $\pair
S{\mathcal U}$ is $\pair{S^n}
{\bigwedge_{i=1}^n(\pi_{in})^{-1}(\mathcal U)}$,
where
$\pi_{in}:S^n\to S$ is the $i^{\text{th}}$ projection.
\end{theorem}

Note that the uniformity has a base consisting of binary
relations of the form
\[
\bigcap_{i=1}^n(\pi^n_i)^{-1}(U_i),
\]
where the $U_i$ are elements of $\mathcal U$.

\subsection{Uniformly
continuous operations}
Let $\pair S{\mathcal U}$ be a uniform space. An $n$-ary
operation on $S$ is \emph{uniformly continuous} if it is
uniformly continuous as a function from
$\pair S{\mathcal U}^n$ to $\pair S{\mathcal U}$.
That is, an $n$-ary operation $\omega$ is uniformly
continuous iff for all $U\in\mathcal U$, there are $U_1$,
$\ldots$, $U_n\in\mathcal U$ such that if
$a_i\mathrel{U_i}b_i$ for all $i$, then $\omega(\vec
a)\mathrel U\omega(\vec b)$.

We
also say that  $\omega$ is \emph{
uniformly continuous at argument $i$} if for all
$U\in\mathcal U$, there is a $\bar U\in\mathcal U$ such
that if $\vec a$,
$\vec b$ are $n$-tuples of $A$ with $a_i\mathrel{\bar
U}b_i$ and $a_j=b_j$ for $j\neq i$, then $\omega(\vec
a)\mathrel U\omega(\vec b\,)$.

\begin{theorem}\label{T:ArgSep}
An $n$-ary operation is uniformly
continuous iff it is uniformly continuous in each
argument.
\end{theorem}

\begin{proof}
Clearly, if an operation is uniformly continuous, it
is uniformly continuous at each argument. To prove the
converse, given
$U$, we use (U5) repeatedly to find that there is a $\hat
U\in\mathcal U$ such that
$\hat U^n\subseteq U$. Then, for each $i$, there is a
$U_i\in\mathcal U$ such that if $a_i\mathrel{
U_i}b_i$ and $a_j=b_j$ for $j\neq i$, then $\omega(\vec
a)\mathrel{\hat U}\omega(\vec b)$. It follows that if
$a_i\mathrel{U_i}b_i$ for all $i$, we have $\omega(\vec
a)\mathrel{\hat U^n}\omega(\vec b)$, whence
$\omega(\vec a)\mathrel U\omega(\vec b)$.
\end{proof}

\subsection{Uniform algebras}
The category of \emph{uniform algebras} in the variety
$\mathbf V$ is the category $\mathbf V[\Unif]$ of $\mathbf
V$-objects in the category $\Unif$. This category can be
defined in two ways, the second of which is more
categorical than the first:

\smallskip

\textbf{Definition 1:} An object of $\mathbf V[\Unif]$
is an algebra in $\mathbf V$, provided with a uniformity
such that all of the basic operations and term
operations are uniformly continuous. An arrow of $\mathbf
V[\Unif]$ is a homomorphism which is uniformly continuous.

\smallskip

\textbf{Definition 2:} An object of $\mathbf V[\Unif]$
is a pair $\pair{\pair S{\mathcal U}}F$, where
$F:\Unif\to\mathbf V$ is a contravariant functor, such
that the contravariant hom functor $\Unif(-,\pair
S{\mathcal U})=U_{\mathbf V}F$, where $U_{\mathbf
V}:\mathbf V\to\Set$ is the underlying set functor on
$\mathbf V$. An arrow
\[
f:\pair{\pair S{\mathcal U}}F\to
\pair{\pair T{\mathcal V}}G
\]
of
$\mathbf V[\Unif]$ is a uniformly continuous function
$f:S\to T$ such that for each uniform space
$\pair Y{\mathcal Y}$, the function
$\Unif(Y,f):\Unif(Y,S)\to\Unif(Y,T)$ is the underlying
function of a homomorphism $\hat f:F(S)\to G(T)$.
\smallskip

We omit the proof that these definitions are equivalent.
See \cite[pp.\ 167--170]{R92} for a fuller
discussion of these two ways of defining $\mathbf
V[\mathbf C]$ for a category $\mathbf C$, and the
advantages of the second definition.

\begin{theorem}\label{T:Compatible} Let $A$ be an algebra,
and let
$\mathcal U\in\OpUnif |A|$.
Then the following are equivalent:
\begin{enumerate}
\item[(1)] $\pair A{\mathcal U}$ is a uniform
algebra.
\item[(2)] Every term operation is uniformly continuous.
\item[(3)] Every basic operation is uniformly continuous.
\item[(4)] Every term operation with arity at least $1$ is
uniformly continuous in the first argument.
\item[(5)] $U^t\in\mathcal U$ for all $U\in\mathcal U$
and for all $n$-ary term operations $t$ with $n>0$.
\item[(6)] $\mathcal U$ is a compatible filter of
relations on $A$.
\end{enumerate}
\end{theorem}

\begin{proof} Clearly (1) $\implies$ (2) $\implies$ (3)
$\implies$ (6), (2) $\implies$ (1), (2) $\implies$ (4)
$\iff$ (5), and (6) $\implies$ (2).

To show that (4) $\implies$ (2), let $t$ be an $n$-ary
term operation. $t$ is trivially uniformly continuous if
$n=0$; suppose $n>0$, and let $1\leq i\leq n$. Let
\[
t'=t\langle\pi_{in},\pi_{1n},\ldots,
\pi_{i-1,n},\pi_{i+1,n},\ldots,\pi_{nn}\rangle
\]
Then $t'$ is uniformly continuous in the first argument
by (4), which implies that $t$ is uniformly continuous in
the $i^{\text{th}}$ argument. Thus, $t$ is uniformly
continuous, by theorem~\ref{T:ArgSep}.
 \end{proof}

\subsection{Compatible uniformities}
If $A$ is an algebra, and $\mathcal U\in\OpUnif |A|$,
then we say that $\mathcal U$ is \emph{compatible} if the
equivalent conditions of theorem~\ref{T:Compatible} are
met. We denote the set of compatible uniformities on $A$
by
$\OpUnif A$.

\begin{theorem} $\OpUnif A$ is a complete meet
subsemilattice of $\OpUnif |A|$.
\end{theorem}

\begin{proof} Follows from theorem \ref{T:CompClosure}.
\end{proof}

\subsection{The underlying algebra functor}
The \emph{underlying algebra} of a uniform algebra $\pair
A{\mathcal U}$ is $A$. The underlying algebra is
functorial in an obvious manner.

\begin{theorem} \label{T:Adjoint} If $A$ is an algebra,
let
$\mathcal U$ be the discrete uniformity on $A$, and let
$\mathcal U'$ be
the indiscrete unifomity. Both $\mathcal U$ and
$\mathcal U'$
are compatible. $\pair{\pair A{\mathcal U}}{1_A}$ is a
universal arrow from $A$ to the underlying algebra
functor from $\mathbf V[\Unif]$ to $\mathbf V$, and
$\pair{\pair A{\mathcal U'}}{1_A}$ is a universal arrow
from
the underlying algebra functor to $A$.
\end{theorem}

\section{Congruentially Uniform
Algebras}\label{S:CongUnif}

Let $A$ be an algebra in the variety $\mathbf
V$. A \emph{congruential uniformity} on $A$ is a
uniformity having a base which is a set of
congruences of $A$. It is easy to see that starting from
any $F\subseteq\Con A$, which is a base for a filter of
congruences of $A$, and viewing it as a base for a filter
of binary relations on
$A$, we obtain a uniformity. We denote the resulting
uniformity by
${\mathcal U}F$, and abbreviate
$\pair A{{\mathcal U}F}$ by $\pair AF$.

Not all uniformities are congruential. For example,
the real numbers are a uniform abelian group, when given
the uniformity having as a base all relations
$U_\epsilon=\{\,\pair xy:|x-y|<\epsilon\,\}$,
but this uniformity is not congruential.

The relationship between filters of congruences and the
congruential uniformities they determine is a one-one
correspondence:

\begin{theorem} \label{T:Base2} A congruence $\alpha$
belongs to the uniformity $\mathcal U$ having as base a
set of congruences
$F$, iff $\alpha$ belongs to the filter of
congruences generated by $F$.\end{theorem}

\begin{proof} $F$ is a base for $\mathcal U$. Thus,
$\alpha\in\mathcal U$ iff $\beta\subseteq\alpha$ for some
$\beta\in F$, iff $\alpha\in F$. \end{proof}

\begin{corollary} The mapping $F\mapsto{\mathcal
U}F$
is a lattice isomorphism from
$\Fil\Con A$ onto a complete meet
subsemilattice of $\OpUnif A$.
\end{corollary}

\begin{corollary} The lattice of congruential
uniformities of an algebra $A$ satisfies every
lattice-theoretic identity satisfied by $\Con A$.
\end{corollary}

\subsection{Congruentialization}
 We denote
by $\CU[\mathbf V]$ the full subcategory of $\mathbf V[\Unif]$
of congruentially uniform algebras. There is an obvious
forgetful functor from $\CU[\mathbf V]$ to $\mathbf V[\Unif]$.

\begin{theorem} If $\pair A{\mathcal U}$ is a uniform
algebra, let $F\subseteq\Con A$ be the filter of
congruences contained in $\mathcal U$. Then $\pair{\pair
AF}{1_A}$ is a universal arrow from $\pair A{\mathcal U}$
to the forgetful functor from
$\CU[\mathbf V]$ to $\mathbf V[\Unif]$. \end{theorem}

\section{The Complete Lattices of
Uniformities}\label{S:Joins}

We have defined three complete lattices of uniformities for
an algebra $A$, namely, the complete lattice $\OpUnif |A|$
of uniformities of the underlying set of $A$, the
complete lattice $\OpUnif A$ of compatible uniformities
of $A$, and the complete lattice of congruential
uniformities, which is isomorphic to $\Fil\Con A$.
So far we have shown that $\Fil\Con A$ is a complete meet
subsemilattice of $\OpUnif A$, and that $\OpUnif A$ is a
complete meet subsemilattice of $\OpUnif |A|$. In this
section, we will examine the join operations and show
that the second of these inclusions is a complete
lattice homomorphism. (For the best result we have
proved regarding the first inclusion, see theorem
\ref{T:FirstMap}.)

\subsection{Dividing sequences and divisible elements}
If $\mathcal U$ is a filter of binary relations on
$S$, we say that a sequence $U_0$, $U_1$,
$\ldots$ of elements of $\mathcal U$ is a
\emph{dividing sequence} in $\mathcal U$ if 
 for each $i$, $U_{i+1}\circ
U_{i+1}\subseteq U_i$. We say that an element
$U\in\mathcal U$ is \emph{divisible} in $\mathcal U$ if
there is a dividing sequence
$U_0$, $U_1$, $\ldots$ in $\mathcal U$ such that $U=U_0$
We denote by
$\Div\mathcal U$ the set of dividing sequences in
$\mathcal U$,
and by $\Div_0\mathcal U$ the set of divisible elements
of $\mathcal U$.

\begin{theorem}
\label{T:Join}
Let $S$ be a set. We have
\begin{enumerate}
\item[(1)] If $\mathcal U$ is a semiuniformity
on $S$, then $\Div_0\mathcal U$ is a uniformity on
$S$, and is the least 
$\mathcal V\in\OpUnif S$ such that
$\mathcal V\subseteq\mathcal U$.
\item[(2)] $\Div_0$, viewed as a function from
$\OpSemiUnif S$ to $\OpUnif S$, preserves arbitrary joins.
\item[(3)] If $\{\,\mathcal U_i\,\}_{i\in I}$ is a
tuple of uniformities on $S$, then the join 
$\bigvee_i\mathcal U_i$, in the lattice $\OpUnif S$, is
$\Div_0(\bigcap_i \mathcal
U_i)$.
\item[(4)] If $\{\,\mathcal U_i\,\}_{i\in I}$ is  a
tuple of uniformities on $S$, then $\Div\bigvee_i\mathcal
U_i=\bigcap_i\Div\mathcal U_i$.
\end{enumerate}
\end{theorem}

\subsection{$\OpUnif A$ and $\Con A$}
\begin{theorem} The mapping
$\alpha\mapsto{\mathcal U}\{\,\alpha\,\}$ is a
lattice isomorphism from $\Con A$ onto a complete join
subsemilattice of $\OpUnif A$.
\end{theorem}

\begin{proof}
Let $\{\,\alpha_i\,\}_{i\in I}$
be a tuple of elements of $\Con A$.
Then $\bigcap_i{\mathcal U}\{\,\alpha_i\,\}$ is
the principal filter generated by the union
$\bigcup_i\alpha_i$. The divisible
elements of this filter all contain
$\bigvee_i\alpha_i$, which is itself divisible. Thus,
$\bigvee_i{\mathcal
U}\{\,\alpha_i\,\}={\mathcal
U}\{\,\bigvee_i\alpha_i\,\}$.
\end{proof}

\subsection{$\OpUnif A$ and $\OpUnif |A|$}
The next theorem parallels the well-known result that
the join of elements of a congruence lattice is the
same as their join in the lattice of equivalence
relations.  Recall the definition of
$U^t$ from section \ref{S:Relations}.

\begin{theorem} Let $A$ be an algebra. If $\mathcal
U_i\in\OpUnif A$ for all $i\in I$, then
\[
\bigvee_i^{\OpUnif
A}\mathcal U_i=\bigvee_i^{\OpUnif|A|}\mathcal U_i.
\]
\end{theorem}

\begin{proof} Denote $\bigvee^{\OpUnif|A|}\mathcal U_i$ by
$\mathcal Y$. We must show that $\mathcal Y$ is compatible,
i.e., that if $t$ is an $n$-ary term for $n>0$, and
$Y\in\mathcal Y$, then
$Y^t\in\mathcal Y$.

Let $Y_0$, $Y_1$,
$\ldots$ be a dividing sequence of elements
of
$\mathcal Y$ with $Y_0=Y$. For each $j$, we have
$Y_j=\bigcup_i U_{ij}$ for some
$U_{ij}\in\mathcal U_i$. We have
$Y_j^t\in\bigcap_i\mathcal U_i$ by lemma~\ref{T:Subset}.
The
$Y^t_j$ form a dividing sequence of elements of
$\bigcap_i\mathcal U_i$, because
\begin{align*}
a\mathrel{Y^t_{j+1}}a'\mathrel{Y^t_{j+1}}a''\\
&\implies
t(a,\vec b\,)\mathrel{Y_{j+1}}t(a',\vec
b\,)\mathrel{Y_{j+1}}t(a'',\vec b\,)
& &\text{for all $\vec b$}
\\
&\implies t(a,\vec b\,)\mathrel{Y_j}t(a'',\vec b\,)
& &\text{for all $\vec b$} \\
&\implies a\mathrel{Y^t_j}a''.
\end{align*}
Thus, $Y^t\in\mathcal Y$.
\end{proof}

\subsection{Relational products and joins in $\OpUnif A$}
It is clear that if $\{\,\mathcal
U_i\,\}_{i\in I}$ is a tuple of uniformities, then
$\bigvee_i\mathcal U_i$ has a base of relations of the
form $\bigcup_iU_i$ with
$U_i\in\mathcal U_i$. Similarly, we have

\begin{lemma} \label{T:Base4} If $\{\,\mathcal
U_i\,\}_{1\leq i\leq n}$ is a finite tuple of
uniformities, then
$\bigvee_i\mathcal U_i$ has a base of relations of the
form $U_1\circ\ldots\circ U_n$ with $U_i\in\mathcal U_i$.
\end{lemma}

Also, we have

\begin{lemma} If $\mathcal U$ is a uniformity, and
$\mathcal V_1$, $\mathcal V_2$ are filters of reflexive
relations such that $\mathcal V_1\leq\mathcal U$,
$\mathcal V_2\leq\mathcal U$, then $\Fg(\mathcal
V_1\circ\mathcal V_2)\leq
\mathcal U$.
\end{lemma}

Now we will give a formula for $\mathcal U\circ
\mathcal V$ which will be useful in section \ref{S:CM}:

\begin{theorem} \label{T:Join2}
Let $\mathcal U$ and
$\mathcal V$ be uniformities on $S$, and define $\mathcal
R_k$ for each cardinal number $k$, inductively as follows:
\begin{enumerate}
\item[$\bullet$] $\mathcal R_0=\mathcal V$;
\item[$\bullet$] $\mathcal R_{k+1}=\mathcal R_k\circ
\mathcal U\circ\mathcal R_k$;
\item[$\bullet$] $\mathcal R_k=\bigcap_{k'<k}\mathcal
R_{k'}$ for $k$ a limit cardinal.
\end{enumerate}
Then $\mathcal U\vee\mathcal V=\bigcap_k\mathcal R_k$.
\end{theorem}

\begin{proof}
It is easy to see that $\mathcal R_k$ is a semiuniformity
for each $k$, the $\mathcal R_k$ are increasing, and that
$\mathcal R_k\leq\mathcal U\vee\mathcal V$. For some $k$,
we must have $\mathcal R_{k+1}=\mathcal R_k$, which
implies that the sequence becomes stationary and that
$\mathcal R_k\circ\mathcal R_k\leq\mathcal R_k$. Thus,
$\mathcal R_k$ is a uniformity and consequently,
$\mathcal U\vee\mathcal V=\mathcal R_k$.
\end{proof}

\section{Permuting Uniformities and Algebras in Mal'cev
Varieties}\label{S:Malcev}

\subsection{Uniformities which permute}
Let $\mathcal U$, $\mathcal V\in\OpUnif S$. We say that
$\mathcal U$ and $\mathcal V$ \emph{permute} if
$\Fg(\mathcal U\circ\mathcal V)=\Fg(\mathcal
V\circ\mathcal U)$. That is,
$\mathcal U$ and
$\mathcal V$ permute iff
\begin{enumerate}
\item[(P1)] for every $U\in\mathcal U$ and
$V\in\mathcal V$, there exist
$\bar U\in\mathcal U$ and $\bar V\in\mathcal V$ such
that $\bar V\circ\bar U\subseteq U\circ V$, and
\item[(P2)] for every $U\in\mathcal U$ and
$V\in\mathcal V$, there exist
$\hat U\in\mathcal U$ and $\hat V\in\mathcal V$ such that
$\hat U\circ\hat V\subseteq V\circ U$.
\end{enumerate}

\begin{theorem} Let $A$ be an algebra, and
$\mathcal U$,
$\mathcal V\in\OpUnif A$. Then
in the lattice $\OpUnif A$, $\mathcal U\vee\mathcal V
=\Fg(\mathcal U\circ\mathcal V)$ iff $\mathcal U$ and
$\mathcal V$ permute.
\end{theorem}

\begin{proof} We always have $\mathcal U\circ\mathcal V\subseteq\mathcal U\vee\mathcal V$ and $\mathcal V\circ\mathcal U\subseteq\mathcal U\vee\mathcal V$. Thus, if $\mathcal U$ and
$\mathcal V$ to not permute, we cannot have $\mathcal U\circ\mathcal V=\mathcal U\vee\mathcal V$.

On the other hand, let us assume that $\mathcal U$ and $\mathcal V$ permute, and we will show that (B2) through (B5) hold
for $\mathcal U\circ\mathcal V$.

(B2): If $U\circ V$, $U'\circ V'\in\mathcal U\circ\mathcal V$,
then let $U''=U\cap U'$ and $V''=V\cap V'$.
We have
$U''\circ V''
\subseteq U\circ V$ and
$U''\circ V''\subseteq U'\circ V'$,
whence
\[
U''\circ V''
\subseteq (U\circ V)\cap(U'\circ V').
\]

(B3) is clearly satisfied.

(B4): If $U\in\mathcal U$ and $V\in\mathcal V$, then
$U^{-1}\in\mathcal U$ and $V^{-1}\in\mathcal V$, and by (P2), there
are $\bar U\in\mathcal U$, $\bar V\in\mathcal V$ such that
\begin{align*}
\bar U\circ\bar V
&\subseteq V^{-1}\circ U^{-1} \\
&=(U\circ V)^{-1}.
\end{align*}

(B5): If $U\in\mathcal U$ and $V\in\mathcal V$, then by
(U2) and (U5) for
$\mathcal U$ and $\mathcal V$, and using (P1), there exist $\hat
U$, $\bar U\in\mathcal U$, and $\hat V$, $\bar V\in\mathcal V$ such
that
\begin{enumerate}
\item[$\bullet$] $\hat U\circ\hat U\subseteq U$;
\item[$\bullet$] $\hat V\circ\hat V\subseteq V$;
\item[$\bullet$] $\bar V\circ\bar U\subseteq \hat
U\circ\hat V$;
\item[$\bullet$] $\bar U\subseteq\hat U$;
\item[$\bullet$] $\bar V\subseteq\hat V$.
\end{enumerate}
Then we have
\begin{align*}
(\bar U\circ\bar V)\circ
(\bar U\circ\bar V)
&=\bar U\circ\bar V\circ\bar U\circ\bar V\\
&\subseteq\bar U\circ\hat U\circ\hat V\circ\bar V\\
&\subseteq U\circ V.
\end{align*}
\end{proof}

\subsection{Algebras with permuting uniformities}
We say that an algebra $A$ has permuting uniformities if
its compatible uniformities permute with each other
pairwise. Note that this implies that $A$ has permuting
congruences.

 \begin{theorem} Let the algebra $A$
have permuting uniformities. Then $\OpUnif A$ is modular.
\end{theorem}

\begin{proof} Let $\mathcal U$, $\mathcal V$, and $\mathcal W\in\OpUnif A$ such that $\mathcal U\leq\mathcal W$. It suffices to
show that
\[
(\mathcal U\vee\mathcal V)\wedge\mathcal W\leq
\mathcal U\vee(\mathcal V\wedge\mathcal W),
\]
as the reverse inequality holds in every lattice.

Let
$U\in\mathcal U$, $V\in\mathcal V$, and $W\in\mathcal W$; it suffices
to show that $U\circ(V\cap W)$ is contained in
the left hand side, because such elements form a base
for the right-hand side. Elements of the form $(\bar
U\circ\bar V)\cap\bar W$, where  $\bar U\in\mathcal U$,
$\bar V\in\mathcal V$, and $\bar W\in\mathcal W$, form a
base for the left hand side. We choose
$\bar U\in\mathcal U$,
$\bar W\in\mathcal W$ such that
\begin{enumerate}
\item[$\bullet$] $\bar W^{-1}\circ\bar W\subseteq W$,
\item[$\bullet$] $\bar U\subseteq U$, and
\item[$\bullet$] $\bar U\subseteq \bar W$ (using the
assumption that $\mathcal U\leq
\mathcal W$).
\end{enumerate}
If $a\mathrel{((\bar U\circ V)\cap\bar W)}b$, then we have
$a\mathrel{\bar U}c\mathrel Vb$ for some $c$, and
$a\mathrel{\bar W}b$. Then because $\bar U\subseteq\bar
W$, we have
$c\mathrel{\bar W^{-1}}a$, which together with
$a\mathrel{\bar W}b$ implies that $c\mathrel Wb$. It
follows that
\[
(\bar U\circ V)\cap\bar W\subseteq
U\circ(V\cap W).
\]
 \end{proof}

\begin{theorem} \label{T:FirstMap} Let the algebra $A$
have permuting uniformities, and let $F_1$, $\ldots$,
$F_n\in\Fil\Con A$. Then $\mathcal
U(\bigvee^{\Fil\Con A}_iF_i)=\bigvee^{\OpUnif
A}_i{\mathcal U}F_i$.\end{theorem}

\begin{proof} Clearly,
\[
\bigvee_i{\mathcal U}F_i
\leq\mathcal U(\bigvee_iF_i);
\]
to prove the reverse inequality, it suffices to show that
if
$\alpha_i\in F_i$ for all $i$, then there is a
$\beta\in\bigvee_iF_i$ such that
$\beta\subseteq\alpha_1\circ\ldots\circ\alpha_n$, because
by theorem~\ref{T:Base4}, relations of the form
$\alpha_1\circ\ldots\circ\alpha_n$ form a base for
$\bigvee_i({\mathcal U}F_i)$.
We simply choose
$\beta=\alpha_1\circ\ldots\circ\alpha_n$: it is a
congruence because $A$ has permuting congruences, and it is
an element of $\bigvee_iF_i$.
\end{proof}

\subsection{Algebras in Mal'cev varieties} A
\emph{Mal'cev variety} is a variety having a
ternary term $p$, called a \emph{Mal'cev term}, satisfying
the identities $p(x,x,y)=p(y,x,x)=y$. If $\alpha$ and
$\beta$ are congruences of an algebra in such a variety, a
famous theorem \cite{Mal} states that
$\alpha\circ\beta=\beta\circ\alpha$.
Similarly,

\begin{theorem} Let $A$ be an algebra in a Mal'cev
variety, and $\mathcal U$, $\mathcal V\in\OpUnif A$.
Then $\mathcal U$ and $\mathcal V$ permute. \end{theorem}

\begin{proof}
Let $p$ be a Mal'cev term for the variety,
and let $U\in\mathcal U$,
$V\in\mathcal V$. There are
$\bar U\in\mathcal U$, $\bar
V\in\mathcal V$ such that
\begin{enumerate}
\item[$\bullet$] if $x\mathrel{\bar U}x'$, then
$p(x,y,z)\mathrel{U}p(x',y,z)$ for all $y$ and $z$, and
\item[$\bullet$] if $z\mathrel{\bar V}z'$, then
$p(x,y,z)\mathrel{V}p(x,y,z')$ for all $x$ and $y$.
 \end{enumerate}
If
$a\mathrel{\bar U}b\mathrel{\bar V}c$, then we have
$a=p(a,b,b)\mathrel{V}p(a,b,c)\mathrel U p(b,b,c)=c$.
Thus,
$\bar U\circ\bar V\subseteq V\circ U$;
this proves (P2). The proof of (P1) is similar.
\end{proof}

\subsection{Algebras in arithmetical varieties}
An \emph{arithmetical variety} is a Mal'cev variety
having a ternary
term $M$ satisfying the identities
\[
M(x,x,y)=M(x,y,x)=M(y,x,x)=x.
\]
The algebras in such a variety all have distributive
congruence lattices.

\begin{theorem} Let $A$ be an algebra in an
arithmetical variety. Then $\OpUnif A$ is distributive.
\end{theorem}

\begin{proof}
Let $\mathcal U$, $\mathcal V$, $\mathcal W\in\OpUnif A$.
It suffices to prove
\[
\mathcal U\wedge(\mathcal V\vee\mathcal W)\leq
(\mathcal U\wedge\mathcal V)\vee(\mathcal U\wedge\mathcal W),
\]
as the reverse inequality holds in any lattice. Using
lemma \ref{T:Base4}, every element of the right-hand side
contains a relation of the form $(U\cap V)\circ(
U\cap W)$, where $U\in\mathcal U$, $V\in\mathcal V$, and
$W\in\mathcal W$. Let $\bar U\in\mathcal U$
($\bar V\in\mathcal V$, $\bar W\in\mathcal W$) be such that
$a\mathrel{\bar U}a'\implies
M(a,c,b)\mathrel{U}M(a',c,b)$ and $b\mathrel{\bar
U}b'\implies M(a,c,b)\mathrel{U}M(a,c,b')$ (respectively,
$c\mathrel{\bar V}c'\implies M(a,c,b)\mathrel VM(a,c',b)$,
$c\mathrel{\bar W}c'\implies M(a,c,b)\mathrel
WM(a,c',b)$), where $M$ is a ternary term satisfying the
above identities.
If $a\mathrel{\bar U\cap(\bar V\circ\bar
W)}b$, this means that $a\mathrel{\bar U}b$ and
$a\mathrel{\bar V}c\mathrel{\bar W}b$ for some $c\in A$.
Then
\[
a\mathrel{U}M(a,c,b)\mathrel{U}b
\]
and
\[
a\mathrel V M(a,c,b)\mathrel Wb,
\]
showing that
\[
a\mathrel{(U\cap V)\circ(U\cap
W)}b;
\]
the desired
inequality follows.
 \end{proof}

\section{Algebras in Congruence-Modular
Varieties}\label{S:CM}

The results of the previous section show  that some 
Mal'cev conditions, on a variety containing an
algebra $A$, influence the structure of the lattice
$\OpUnif A$ in the same way they influence that of $\Con
A$. 
Now, a theorem
\cite{Day} states that a variety
$\mathbf V$ is congruence-modular iff there is a
finite sequence $m_0$, $\ldots$, $m_n$ of quaternary
terms, called \emph{Day terms}, satisfying the following
identities:
\begin{enumerate}
\item[(D1)] For all $i$, $m_i(x,y,y,x)=x$;
\item[(D2)] $m_0(x,y,z,w)=x$;
\item[(D3)] $m_n(x,y,z,w)=w$;
\item[(D4)] for even $i<n$,
$m_i(x,x,y,y)=m_{i+1}(x,x,y,y)$; and
\item[(D5)] for odd $i<n$,
$m_i(x,y,y,z)=m_{i+1}(x,y,y,z)$.
\end{enumerate}
We are led to ask whether this Mal'cev condition on the
variety $\mathbf V$ forces $\OpUnif A$ to be modular, for
$A\in\mathbf V$.
In this section, we will prove a partial result in this
direction.

For this discussion, $\mathbf V$ will be a
congruence-modular variety, $m_i$, $i=0$, $\ldots$, $n$
will be a sequence of Day terms for $\mathbf V$, and $A$
will be an algebra in $\mathbf V$. Our first goal will be
to prove a generalization of the Shifting Lemma
\cite{Gumm}.

\begin{lemma} \label{T:Lemma1}
Let $\mathcal X\in\OpUnif A$. Given
$X\in\mathcal X$, there exists $\bar X\in\mathcal X$
such that if
$a$, $b$, $c$, $d\in A$ with $b\mathrel{\bar X}d$ and
$m_i(a,a,c,c)\mathrel{\bar X}m_i(a,b,d,c)$ for all $i$,
then
$a\mathrel Xc$.
\end{lemma}

\begin{proof}
For the given $X$, let $\hat X\in\mathcal X$ be such that
$\hat X^{2n}\subseteq X$. Let $\bar X=\bigcap_i\bar X_i$
where each $\bar X_i$ is symmetric and such that
$m_i(\bar X_i)\subseteq\hat X$. Let $a$, $b$, $c$, $d$ be
given such that $b\mathrel{\bar X}d$ and
$m_i(a,a,c,c)\mathrel{\bar X}m_i(a,b,d,c)$ for all $i$.
Define $u_i=m_i(a,b,d,c)$ and $v_i=m_i(a,a,c,c)$ for all
$i$. Since $b\mathrel{\bar X}d$ and $v_i\mathrel{\bar
X}u_i$ for all
$i$, we have for even $i$, $u_i\mathrel{\hat
X}v_i=v_{i+1}\mathrel{\hat X}u_{i+1}$, while for odd $i$,
$u_i\mathrel{\hat
X}m_i(a,b,b,c)=m_{i+1}(a,b,b,c)\mathrel{\hat X}u_{i+1}$.
If follows that for all $i<n$, we have $u_i\mathrel{(\hat
X\circ\hat X)}u_{i+1}$. Thus, $a=u_0\mathrel Xu_n=c$.
\end{proof}

\begin{lemma}\label{T:Lemma2} For every binary relation
$U$ on
$A$, define $M(U)=\bigcup_im_i(U)$. Then if $U$ and
$V$ are binary relations on $A$, we have $M(U\circ
V)\subseteq M(U)\circ M(V)$.
\end{lemma}

\begin{proof} Follows from lemma~\ref{T:Composition}.
\end{proof}

\begin{lemma}\label{T:Lemma3}
Let $\mathcal T_0\in\OpSemiUnif A$ and
$\mathcal T_1$,
$\mathcal X\in\OpUnif A$ be such that
$\mathcal T_0\wedge\mathcal T_1\leq\mathcal X$. Given
$X\in\mathcal X$, there exist $T_0$, $T_1$,  $\bar X$ in
$\mathcal T_0$,
$\mathcal T_1$, $\mathcal X$ respectively, such that $T_0$
and $T_1$ are symmetric and
$a\mathrel{T_0}b\mathrel{T_1}d\mathrel{T_0}c\mathrel{T_1}
a$
and $b\mathrel{\bar X}d$ imply $a\mathrel Xc$.
\end{lemma}

\begin{proof} Let $\bar X\in\mathcal X$
 satisfy the conclusion
of lemma~\ref{T:Lemma1}. Let $\hat T_0\in\mathcal T_0$,
$\hat T_1\in\mathcal T_1$ be symmetric and such that
$\hat T_0\cap\hat T_1\subseteq\bar X$. Let
$T_0\in\mathcal T_0$ be symmetric and such that
$M(T_0)\subseteq\hat T_0$. Using lemma~\ref{T:Lemma2}, let
$T_1\in\mathcal T_1$ be symmetric and such that
$M(T_1)\circ M(T_1)\subseteq\hat T_1$. Then
$m_i(a,a,c,c)\mathrel{\hat T_0}m_i(a,b,d,c)$ for all $i$,
and
$m_i(a,a,c,c)\mathrel{M(T_1)}m_i(a,a,a,a)=a=m_i(a,b,b,a)
\mathrel{M(T_1)}m_i(a,b,d,c)$ for all $i$, whence
$m_i(a,a,c,c)\mathrel{\bar X}m_i(a,b,d,c)$
for all $i$. Then if
$b\mathrel{\bar X}d$, it follows from lemma~\ref{T:Lemma1}
that
$a\mathrel Xc$.
\end{proof}

\begin{lemma}[Shifting Lemma] \label{T:Lemma4}
Let $\mathcal R\in\OpSemiUnif A$ and
$\mathcal W$,
$\mathcal X\in\OpUnif A$ be such that
$\mathcal R\wedge\mathcal W\subseteq\mathcal X$. Then
$(\mathcal R\circ(\mathcal W\wedge\mathcal X)\circ\mathcal
R)\wedge\mathcal W\leq\mathcal X$.
\end{lemma}

\begin{proof} Let $X\in\mathcal X$. Write $\mathcal U$
for $\mathcal W\wedge\mathcal X$. By
lemma~\ref{T:Lemma3}, there exist
$R\in\mathcal R$, $W\in\mathcal W$, and $\bar
X\in\mathcal X$ such that $W$ is symmetric and $a\mathrel
R b\mathrel W d\mathrel R c\mathrel Wa$ and
$b\mathrel{\bar X}d$ imply
$a\mathrel Xc$. Let $U=W\cap\bar X$.  Then
\[a\mathrel{(R\circ U\circ R)\cap W}c\]
implies that $a\mathrel{W}c$ and that there exist
$b$, $d\in A$ such that
$a\mathrel{R}b\mathrel Ud\mathrel{R}c$, which implies
that $a\mathrel R b\mathrel Wd\mathrel R
c\mathrel Wa$ and that $b\mathrel{\bar X}d$. We then have
$a\mathrel Xc$. Thus, $(R\circ U\circ R)\cap
W\subseteq X$, implying that $(\mathcal R\circ\mathcal
U\circ\mathcal R)\wedge \mathcal W\leq\mathcal X$.
\end{proof}

Now, although we have not proved that $\OpUnif A$ is
modular, we have the following partial result:

\begin{theorem} Let $\mathcal U$, $\mathcal V$,
$\mathcal W\in\OpUnif A$ where $\mathcal
W={\mathcal U}\{\,W\,\}$ for a congruence $W$, and such
that $\mathcal U\leq\mathcal W$. Then
\[\mathcal U\vee(\mathcal V\wedge\mathcal W)
=(\mathcal U\vee\mathcal V)\wedge\mathcal W.
\]
\end{theorem}

\begin{proof}
It is trivial that $\mathcal U\vee(\mathcal
V\wedge\mathcal W)\leq (\mathcal U\vee\mathcal
V)\wedge\mathcal W$, as this inequality holds in every
lattice.

Define $\mathcal R_k$ for every cardinal $k$ as $\mathcal
R_0=\mathcal V$,
$\mathcal R_{k+1}=\mathcal R_k\circ \mathcal
U\circ\mathcal R_k$, and for limit cardinals $k$,
$\mathcal R_k=\bigcap_{k'<k}\mathcal R_{k'}$. By
lemma~\ref{T:Lemma4} and because $\mathcal W$ is generated
by a congruence $W$, we have $\mathcal R_k\wedge\mathcal
W\leq\mathcal\mathcal U\vee(\mathcal V\wedge\mathcal W)$
for all
$k$. By theorem \ref{T:Join2}, we have $\mathcal
U\vee\mathcal V=\mathcal R_k$ for some cardinal $k$.
Thus,
$(\mathcal U\vee\mathcal V)\wedge\mathcal W\leq
\mathcal U\vee(\mathcal V\wedge\mathcal W)$.
\end{proof}

\section{Hausdorff Uniform Algebras}\label{S:Hausdorff}

\subsection{Hausdorff uniform spaces}
Recall the following
well-known proposition:

\begin{proposition} Let $\pair S{\mathcal U}$ be a uniform
space, with associated topology $\mathbf T$. Then the
following are equivalent:
\begin{enumerate}
\item[(1)] $\mathbf T$ is $T_1$;
\item[(2)] $\mathbf T$ is Hausdorff;
\item[(3)] $\bigcap\mathcal U=\Delta_S$.
\end{enumerate}
\qed\end{proposition}

If these equivalent properties are satisfied, then the
uniform space $\pair S{\mathcal U}$ is said to be \emph{
Hausdorff}. We denote the full subcategory of $\Unif$,
with objects the Hausdorff uniform algebras, by
$\HUnif$.

 \subsection{Hausdorffization}
It is not hard to impose the Hausdorff property on a
uniform space:

 \begin{proposition} Let $\pair S{\mathcal U}$ be
a uniform space. Let $\psi=\bigcap\mathcal U$. We have
\begin{enumerate}
\item[(1)] $\psi$ is an
equivalence relation;
\item[(2)] the set $T=S/\psi$, provided
with the direct image uniformity $\mathcal
V=(\nat\psi)_*(\mathcal U)$, is Hausdorff; and
\item[(3)] $\pair{\pair T{\mathcal V}}{\nat\psi}$ is a
universal arrow from $S$ to the forgetful functor from
$\HUnif$ to $\Unif$.
\end{enumerate}
\end{proposition}

\begin{proof} All elements of $\mathcal U$ are reflexive, and
$\mathcal U$ has a base of symmetric relations, so
$\bigcap\mathcal U$ is reflexive and symmetric. That
$\bigcap\mathcal U$ is transitive follows from the inclusion
\[
\bigcap\mathcal U\circ\bigcap\mathcal U\subseteq\bigcap\mathcal U,
\]
which follows from property (U5) of $\mathcal U$.

That $\pair{S/\psi}{\mathcal V}$ is Hausdorff
then follows from theorem~\ref{T:DirectImage} on the form
of the direct image uniformity.

If $f:\pair S{\mathcal U}\to\pair{S'}{\mathcal
U'}$ is uniformly continuous, and
$\psi'=\bigcap\mathcal U'$, then
$f^{-1}(\psi')\geq\psi$. Thus, $f$ must factor
through $\nat\psi$, implying (3).
\end{proof}

\subsection{Hausdorffization and uniform
algebras}
 A uniform algebra $\pair A{\mathcal U}$ is
\emph{Hausdorff} iff it is an object of $\mathbf
V[\HUnif]$. (Meaning simply that the underlying
uniform space of
$\pair A{\mathcal U}$ is an object of $\HUnif$.)

\begin{theorem}
Let $\pair A{\mathcal U}$ be an object of $\mathbf V[\Unif]$.
We have
\begin{enumerate}
\item[(1)] The equivalence relation
$\psi=\bigcap\mathcal U$ is a congruence of $A$;
\item[(2)] $A/\psi$, provided with the direct image
uniformity
$\mathcal V=(\nat\psi)_*(\mathcal U)$, is a uniform
algebra; and
\item[(3)] $\pair{\pair{A/\psi}{\mathcal V}}{\nat\psi}$
is a universal arrow from $A$ to the forgetful functor
from
$\mathbf V[\HUnif]$ to $\mathbf V[\Unif]$.
\end{enumerate}
\end{theorem}

\begin{proof}
Let $t$ be an $(n+1)$-ary term operation of the variety,
let $\vec b\in A^n$, and let $a$, $a'\in A$ be such that
$a\mathrel
\psi a'$. Since $a\mathrel\psi a'$, we have
$a\mathrel{U^t}a'$ for every
$U\in\mathcal U$. This implies that
$t(a,\vec b\,)\mathrel U t(a',\vec b\,)$. As this holds
for all
$U\in\mathcal U$,
$t(a,\vec b)\mathrel\psi t(a',\vec b)$. Thus, $\psi$ is a
congruence.

To prove $\pair{A/\psi}{\mathcal V}$ is a uniform algebra,
let $\omega$ be an $n$-ary basic operation
for $n>0$, and let $U\in\mathcal U$. Let $U_1$,
$\ldots$, $U_n\in\mathcal U$ be such that $\omega(\vec
a)\mathrel U\omega(\vec b\,)$ if $a_i\mathrel{U_i}b_i$
for all $i$. Then $\omega^{A/\psi}(\vec
c\,)\mathrel{V_U}\omega^{A/\psi}(\vec d\,)$ if
$c_i\mathrel{V_{U_i}}d_i$ for all $i$, where
$V_U=(\nat\psi)(U)$. Thus,
$\omega^{A/\psi}$ is uniformly continuous, and
$\pair{A/\psi}{\mathcal V}$ is a uniform algebra.
\end{proof}

As Hausdorffization gives rise to universal arrows, it is
functorial, and we denote the Hausdorffization
functor by $H$.
Clearly, $H$ takes the subcategory $\CU[\mathbf V]$ of
$\mathbf V[\Unif]$ into $\CU[\mathbf V]$; we
denote the subcategory of $\mathbf V[\Unif]$ of
Hausdorff, congruentially uniform algebras by $\HCU[\mathbf
V]$.

\section{Complete Uniform Algebras}\label{S:Complete}

\subsection{Convergence}
Let $\pair S{\mathcal U}$ be a uniform space. A \emph{net}
in $\pair S{\mathcal U}$ is a function from a directed set
$\mathbf D$ into $S$.
We say
that a net $N$ \emph{converges}
to $s\in S$ (with respect to $\mathcal U$) if for each
$U\in\mathcal U$, $\pair s{N(d)}\in U$ for large enough $d$,
i.e., if there is a $\bar d\in\mathbf D$ such that $d\geq
\bar d$ implies $s\mathrel UN(d)$. In this case, we say
$N\rightarrow
s$.

A net $N:\mathbf D\to S$ is a \emph{constant net} (at
$s\in S$) if $N(d)=s$ for all $d\in\mathbf D$.

\begin{proposition}
Let $\pair S{\mathcal U}$ be a uniform space. We have
\begin{enumerate}
\item[(1)] If $N$ is a constant net
at
$s$, then $N\rightarrow s$;

\item[(2)] if $N$ is a net in $S$ such that $N\rightarrow
s$, $\pair T{\mathcal V}$ is another uniform space,
and $f:S\to T$ is a uniformly continuous function, then
$f\circ N\rightarrow f(s)$; and
\item[(3)] if $S'\subseteq S$ is a
dense subset, and $s\in S$, then there is a net $N$ with
values in $S'$ such that $N\rightarrow s$.
 \end{enumerate}
 Now, let $N$ be a net in $S$, and let
$\psi=\bigcap\mathcal U$. We have
\begin{enumerate}
\item[(4)] If $N\rightarrow s$ and $s'\in S$, then
$N\rightarrow s'$ iff $s\mathrel\psi s'$.
\end{enumerate}
\end{proposition}

\begin{corollary}
\label{T:Unique}
If $\pair S{\mathcal U}$ is Hausdorff, then a net
converges to at most one element of $S$, and conversely,
if no net converges to more than one element of
$S$, then $\pair S{\mathcal U}$ is Hausdorff.
\end{corollary}

\begin{proof} Follows from part (4) of the theorem.
\end{proof}

\subsection{Cauchy nets}
Let $\pair S{\mathcal U}$ be a uniform space. A net $N:\mathbf
D\to S$ is a \emph{Cauchy net} (with respect to $\mathcal
U$) if for each $U\in\mathcal U$, $N(d)\mathrel UN(d')$
for large enough $d$ and $d'$, i.e., if there is a $\bar
d\in\mathbf D$ such that $\pair {N(d)}{N(d')}\in U$ whenever
$d\geq\bar d$ and $d'\geq\bar d$. 

\begin{proposition}
\label{T:Cauchy}
We have
\begin{enumerate}
\item[(1)] If $\pair S{\mathcal U}$ is a uniform space,
$s\in S$, and $N:\mathbf D\to S$ is a net such that
$N\rightarrow s$, then $N$ is a Cauchy net.
\item[(2)] If $\pair S{\mathcal U}$ and $\pair T{\mathcal
V}$ are uniform spaces, $f:S\to T$ is a uniformly
continuous function, and $N$ is a Cauchy net in $\pair
S{\mathcal U}$, then $f\circ N$ is a Cauchy net in $\pair
T{\mathcal V}$.
 \item[(3)]
Let $S$ be a set, $\pair T{\mathcal V}$ a uniform space, and $f:S\to T$ a function. Then
$f\circ N$ is a Cauchy net with respect to $\mathcal V$ iff
$N$ is a Cauchy net with respect to $f^{-1}(\mathcal V)$.
\item[(4)] If $S$ is a set, $\mathcal U$, $\mathcal
V\in\OpUnif S$ are such that $\mathcal U\leq\mathcal V$,
and $N$ is a Cauchy net in $S$ with respect to $\mathcal
U$, then $S$ is also Cauchy with respect to $\mathcal V$.
\end{enumerate}
 \end{proposition}

 If every Cauchy net in a uniform space $\pair
S{\mathcal U}$ converges to an element of $S$, then we say that
$\pair S{\mathcal U}$ is \emph{complete}. Complete,
Hausdorff uniform spaces are particularly important to us
and we denote the full subcategory of $\Unif$ of such
spaces by
$\CHUnif$.

\subsection{Equivalence of Cauchy nets}
We say that
Cauchy nets $N:\mathbf D\to S$ and $N':\mathbf D'\to S$
in $\pair S{\mathcal U}$ are \emph{equivalent}, or $N\sim
N'$, if for each $U\in\mathcal U$, $N(d)\mathrel UN'(d')$ for
large enough $d$ and $d'$.

\begin{proposition}
\label{T:Equivalence}
Let $\pair S{\mathcal U}$ be a
uniform space. We have
\begin{enumerate}
\item[(1)] Equivalence of Cauchy nets in $\pair
S{\mathcal U}$ is an equivalence relation.
\end{enumerate}
Now, let $N$, $N'$ be Cauchy nets in $\pair S{\mathcal U}$. We
have
\begin{enumerate}

\item[(2)] If $N\rightarrow s$, and
$N'\sim N$, then $N'\rightarrow s$; and

\item[(3)] if $N\rightarrow s$, $N'\rightarrow s'$, and
$s\mathrel\psi s'$ (in particular, if $s=s'$) where $\psi
=
\bigcap\mathcal U$, then
$N\sim N'$.

\end{enumerate}
Now let $\pair T{\mathcal V}$ be another uniform space, and
$f:S\to T$ a uniformly continuous function. We have
\begin{enumerate}
\item[(4)] If $N\sim N'$, then $f\circ N\sim
f\circ N'$.

\end{enumerate}
\end{proposition}

\subsection{Hausdorff completion of uniform
spaces} If $\pair S{\mathcal U}$ is a uniform
space, let $S/\mathcal U$ be the set of of equivalence
classes of Cauchy nets with respect to $\mathcal U$.

\begin{proposition} $S/\mathcal U$ is a small set.
\end{proposition}

\begin{proof} Every Cauchy net in $\pair S{\mathcal U}$
can be shown equivalent to a net $N:\mathbf D\to S$ where
the directed set $\mathbf D$ is $\mathcal U$, ordered by
reverse inclusion.
\end{proof}

Let $\pair S{\mathcal U}$ be a uniform space, and $V$ a
binary relation on $S$. Let $R(\mathcal U,V)$ be the
binary relation on $S/\mathcal U$ defined by
\[
k\mathrel{R(\mathcal U,V)}k'\iff\exists N\in k,\,N'\in
k'
\text{ such that }N(d)\mathrel VN'(d')\text{ for large
enough $d$, $d'$}
\]
and let $\bar R(\mathcal U,V)$ be the binary relation on
$S/\mathcal U$ defined by
\[
k\mathrel{\bar R(\mathcal U,V)}k'\iff\forall N\in
k,\,N'\in k' ,\,N(d)\mathrel VN'(d')\text{ for large
enough $d$, $d'$\quad}
\]

\begin{lemma} Let $\pair S{\mathcal U}$ be a uniform
space,
and let $V$ be a binary relation on $S$. We have
\begin{enumerate}
\item[(1)] $R(\mathcal U,V)^{-1}=R(\mathcal U,V^{-1})$;
\item[(2)] $\bar R(\mathcal U,V)^{-1}=\bar R(\mathcal U,
V^{-1})$; and
\item[(3)] $\bar R(\mathcal U,V)\subseteq R(\mathcal
U,V)$.
\end{enumerate}
Now, let $V'$ be another binary relation on $S$, such that
$V\subseteq V'$. We have
\begin{enumerate}
\item[(4)] $R(\mathcal U,V)\subseteq R(\mathcal U,V')$;
and
\item[(5)] $\bar R(\mathcal U,V)\subseteq\bar R(\mathcal
U,V')$.
\end{enumerate}
 Now, let $\tilde V\in\mathcal U$ be such that $\tilde
V^3\subseteq V$. We have
\begin{enumerate}
\item[(6)] $R(\mathcal U,\tilde V)\subseteq\bar
R(\mathcal U,V)$; and
\item[(7)] $\bar R(\mathcal U,\tilde V)^2\subseteq \bar
R(\mathcal U,V)$.
\end{enumerate}
\end{lemma}

\begin{theorem}
\label{T:Base3}
Let $\pair S{\mathcal U}$ be a uniform space,
$\mathcal V$ a uniformity on $S$ such that $\mathcal
U\leq\mathcal V$, and $\mathcal B$ a base for $\mathcal
V$. We have
\begin{enumerate}
\item[(1)] The set $\mathcal B_1=\{\,R(\mathcal U,V)\mid
V\in\mathcal B\,\}$ is a base for a uniformity $\mathcal
V/\mathcal U$ of $S/\mathcal U$;
\item[(2)] The set $\mathcal B_2=\{\,\bar R(\mathcal
U,V)\mid V\in\mathcal B\,\}$ is also a base for
$\mathcal V/\mathcal U$; and
\item[(3)] $\mathcal V/\mathcal U$ is independent of the
base $\mathcal B$.
\end{enumerate}
\end{theorem}

\begin{proof}
$\mathcal B_1$ and $\mathcal B_2$ are bases for filters
of relations on $S/\mathcal U$, and the filters they
generate are independent of the chosen base of $\mathcal
V$, by parts (4) and (5) of the lemma.

We must verify (B3), (B4), and (B5) for $\mathcal B_1$
and $\mathcal B_2$. Let $V\in\mathcal B$ and $k\in
S/\mathcal U$. For all $N$, $N'\in k$, we have
$N(d)\mathrel VN'(d')$ for large enough $d$ and $d'$,
because $V\in\mathcal U$. Thus $k\mathrel{\bar R(\mathcal
U,V)}k$, verifying (B3) for $\mathcal B_2$ and, by part
(3) of the lemma, for $\mathcal B_1$.

By parts (1) and (2) of the lemma, and (B4) for $\mathcal
B$, we have (B4) for $\mathcal B_1$ and $\mathcal B_2$.

$\mathcal B_1$ and $\mathcal B_2$ generate the same
filter by parts (3) and (6) of the lemma, while $\mathcal
B_2$ (and therefore also $\mathcal B_1$) satisfies (B5) by
part (7).
\end{proof}

 We define $C\pair S{\mathcal U}$, the \emph{
Hausdorff completion} of $\pair S{\mathcal U}$, to be
$\pair{S/\mathcal U}{\mathcal U/\mathcal U}$.
Let
 $\eta_{\mathcal U}:S\to S/\mathcal U$ be the
function taking an element $s$
to the equivalence class of Cauchy nets converging to $s$.
Most parts of the following theorem are well known:

\begin{theorem}
\label{T:Completion}
Let $\pair S{\mathcal U}$ be a uniform
space. We
have
\begin{enumerate}

\item[(1)] $C\pair S{\mathcal U}$ is Hausdorff and
complete;

\item[(2)] if $\mathcal V\in\OpUnif S$ is such that
$\mathcal U\leq\mathcal V$, then $\eta^{-1}_{\mathcal U}(\mathcal V/\mathcal U)=\mathcal V$;

\item[(3)] $\eta_{\mathcal U}$ is a uniformly continuous
function from $\pair S{\mathcal U}$ to $C\pair S{\mathcal
U}$;

\item[(4)] the
image of $\eta_{\mathcal U}$ is dense in $C\pair S{\mathcal U}$;

\item[(5)] if $\mathcal W\in\OpUnif S/\mathcal U$ is such
that
$\mathcal U/\mathcal U\leq\mathcal W$, then $\mathcal W=\eta^{-1}_{\mathcal U}(\mathcal W)/\mathcal U$;

\item[(6)]
$\pair{C\pair S{\mathcal U}}{\eta_{\mathcal U}}$ is a universal arrow from
$\pair S{\mathcal U}$ to the forgetful functor from $\CHUnif$
to $\Unif$;

\item[(7)]
if $S$ is a set, $\pair T{\mathcal V}$ is a uniform
space with $\pair T{\mathcal V}$ complete and Hausdorff, and
$f:S\to T$ is a function such that
the image of $f$ is dense in $T$, then $\pair T{\mathcal V}$
is isomorphic to $C\pair S{f^{-1}(\mathcal V)}$;

\item[(8)] if $\mathcal V\in\OpUnif S$ is such that
$\mathcal U\leq\mathcal V$, then there is an isomorphism
$g_{\mathcal U,\mathcal V}:C\pair S{\mathcal V}\to
C\pair{S/\mathcal U}{\mathcal V/\mathcal U}$ such that
$\eta_{\mathcal V/\mathcal U}\eta_{\mathcal
U}=g_{\mathcal U,\mathcal V}\eta_{\mathcal V}$; and

\item[(9)] if $\mathcal V$, $\mathcal W\in\OpUnif S$ are
such that
$\mathcal U\leq\mathcal V\leq\mathcal W$, then $g_{\mathcal U,\mathcal V}^{-1}((\mathcal W/\mathcal U)/(\mathcal V/\mathcal U))=\mathcal W/\mathcal V$.
\end{enumerate}
\end{theorem}

\begin{proof} (1): Let $k$, $k'\in S/\mathcal U$ be such 
that
$k\mathrel{\bar R(\mathcal U,U)}k'$ for every
$U\in\mathcal U$, and let
$N\in k$ and $N'\in k'$. For every $U\in\mathcal U$,
$N(d)\mathrel UN'(d')$ for large enough $d$ and $d'$.
Since
$U$ is arbitrary we have
$N\sim N'$, which implies $k=k'$. Thus, 
$\bigcap(\mathcal U/\mathcal U)=\Delta_{S/\mathcal
U}$, implying that
$C\pair S{\mathcal U}$ is Hausdorff.

To prove $C\pair S{\mathcal U}$ is complete, let $K:\mathbf
D\to S/\mathcal U$ be a net in $S/\mathcal U$, Cauchy with
respect to $\mathcal U/\mathcal U$. For each $d\in\mathbf D$, let
$N_d:\mathbf D_d\rightarrow S$ be a representative of
$K(d)$.
For any $U\in\mathcal U$,
for large enough $d$, $d'\in\mathbf D$, $N_d(\bar
d)\mathrel UN_{d'}(\bar d')$ for large enough $\bar
d\in\mathbf D_d$ and $\bar d'\in\mathbf D_{d'}$.
For, we have $K(d)\mathrel{\bar R(\mathcal U,U)}K(d')$ for
large enough $d$ and $d'$.

Consider
$\mathcal U$ as a directed set with ordering given by
reverse inclusion, and for each $U\in\mathcal U$, let
$d(U)\in\mathbf D$, $\bar d(U)\in\mathbf D_{d(U)}$ be
such that $N_{d(U)}(\bar d(U))\mathrel UN_{d'}(\bar d')$
for
$d'\geq d(U)$ and large enough $\bar d'$, and let
$N(U)=N_{d(U)}(\bar d(U))\in S$.
We claim that $N:\mathcal U\to S$ is a Cauchy net in $S$
with respect to $\mathcal U$. To prove this, let
$U\in\mathcal U$. Let $\bar U\in\mathcal U$ be such that
$\bar U\circ\bar U\subseteq U$. Then if $U_1$,
$U_2\subseteq\bar U$, we have
\[
N(U_1)=N_{d(U_1)}(\bar d(U_1))\mathrel{\bar U}
N_{\hat d}(\tilde d)\mathrel{\bar U}N_{d(U_2)}(\bar
d(U_2))=N(U_2)
\]
for some $\hat d$ such that $\hat d\geq d(U_1)$ and $\hat
d\geq d(U_2)$, and some $\tilde d\in\mathbf D_{\hat d}$,
implying that $N(U_1)\mathrel UN(U_2)$. It follows that
$N$ is a Cauchy net.

Let $k$ be the class of Cauchy nets equivalent to
$N$. We will show $K\rightarrow k$, for which it suffices
to show that for all $U\in\mathcal U$, 
$K(d)\mathrel{R(\mathcal U,U)}k$ for large enough $d$.
Again, let $\bar U\in\mathcal U$ be such that $\bar
U\circ\bar U\subseteq U$. If $d\geq d(\bar U)$ and $\bar
d$ is large enough, then
\[
N_d(\bar d)\mathrel{\bar U}N_{\hat d}(\tilde
d)\mathrel{\bar U}N_{d(U')}(\bar d(U'))=N(U'),
\]
whenever $\hat d\geq d(\bar U)$, $U'\subseteq\bar U$, and
$\tilde d$ is large enough. That is,
$K(d)\mathrel{R(\mathcal U,U)}k$ for large enough $d$.
Thus, $S/\mathcal U$ is complete with respect to
$\mathcal U/\mathcal U$.

(2): Given $V\in\mathcal V$, let $\bar V\in\mathcal V$
be such that $\bar V\circ\bar V\circ\bar V\subseteq V$.
Now suppose that $\eta_{\mathcal U}(s)\mathrel{R(\mathcal
U,\bar
V)}\eta_{\mathcal U}(s')$, i.e., that there are sequences
$N$ and $N'$, converging to $s$ and $s'$ respectively,
such that $N(d)\mathrel{\bar V}N'(d')$ for large enough
$d$ and
$d'$. Then since $N$ and $N'$ converge to $s$ and $s'$
respectively, we have $s\mathrel{\bar V}N(d)\mathrel{\bar
V}N'(d')\mathrel{\bar V}s'$ for large enough $d$ and
$d'$, implying $s\mathrel Vs'$. Thus, $\eta^{-1}_{\mathcal U}(\mathcal V/\mathcal U)\leq\mathcal V$.
On the other hand, suppose $s\mathrel{\bar V}s'$. Then if
$N$, $N'$ converge to $s$, $s'$ respectively, we have
$N(d)\mathrel{\bar V}s\mathrel{\bar V}s'\mathrel{\bar
V}N'(d')$ for large enough $d$ and $d'$, showing that
$s\mathrel{\eta^{-1}_{\mathcal U}(\bar R(\mathcal
U,V))}s'$. Thus, $\mathcal V\leq\eta_{\mathcal
U}^{-1}(\mathcal V/\mathcal U)$.

(3): By (2), $\eta^{-1}_{\mathcal U}(\mathcal U/\mathcal U)=\mathcal U$,
implying $\eta_U$ is uniformly continuous.

(4): Given $k\in C\pair S{\mathcal U}$, and $U\in\mathcal
U$, let $N\in k$ where $N:\mathbf D\to S$. Let
$ d\in\mathbf D$ be such that $N(d')\mathrel U N(d'')$ for
$d'$, $d''\geq d$, let $s=N(d)$, and let $N_s$ be a
constant net at $s$. Then $N(d')\mathrel UN_s(d'')$ for
large enough 
$d'$ and $d''$, or $k\mathrel{R(\mathcal U,U)}\eta_U(s)$. 
Thus,
$\eta_{\mathcal U}(S)$ is dense in $C\pair S{\mathcal U}$.

(5): By (2),  $\eta^{-1}_{\mathcal U}(\eta^{-1}_{\mathcal
U}(\mathcal W)/\mathcal U)=\eta^{-1}_{\mathcal
U}(\mathcal W)$. Thus, by (4) and corollary
\ref{T:Dense2}, 
$\eta^{-1}_{\mathcal U}(\mathcal W)/\mathcal U=\mathcal
W$.

(6): We
must show that if $f:\pair S{\mathcal U}\to\pair 
T{\mathcal V}$ is uniformly continuous, where $\pair
T{\mathcal V}$ is complete and Hausdorff, then there is a
unique uniformly continuous $g:\pair{S/\mathcal
U}{\mathcal U/\mathcal U}\to\pair T{\mathcal V}$ such that
$f=g\eta_{\mathcal U}$.

If $N\in k\in S/\mathcal U$, then $f\circ N$ is a Cauchy
net in $\pair T{\mathcal V}$, and we define $g(k)$ to be
the limit of $f\circ N$ in $\pair T{\mathcal V}$. This is
well-defined by corollary \ref{T:Unique} and theorem
\ref{T:Equivalence}. Clearly
$f=g\eta_{\mathcal U}$.

$g$ is unique, because $\eta_{\mathcal U}$, having
dense image, is an epimorphism.

To show $g$ is uniformly continuous, let
$V\in\mathcal V$ and let $U\in\mathcal U$,
$\bar V\in\mathcal V$ be such that $\bar V\circ f(U)\circ
\bar V\subseteq V$. Let $k$ and
$k'$ be elements of $S/\mathcal U$ with
$k\mathrel{\bar R(\mathcal U,U)}k'$; then for some
$N\in k$, $N'\in k'$, $d$,
$d'$ we have $g(k)\mathrel{\bar V}N(d)\mathrel U N'(
d')\mathrel{\bar V}g(k')$, whence
$g(k)\mathrel{V} g(k')$. Thus,
$g(\bar R(\mathcal U,U))\subseteq V$.

(7): Suppose $S$, $\pair T{\mathcal V}$, and $f$ are given as stated. Let $\mathcal W=f^{-1}(\mathcal V)$. By (6), there is a unique, uniformly
continuous function $g:\pair{S/\mathcal W}{\mathcal W/\mathcal W}\to\pair T{\mathcal V}$ such that
$f=g\eta_{\mathcal W}$. We must show that $g$ is a
uniform isomorphism.

If $t\in T$, then since $f(S)$ is dense in $T$, there is
a net $N$ in $S$ such that $f\circ N\rightarrow t$. $N$
is Cauchy with respect to $\mathcal W$ by theorem
\ref{T:Cauchy}(3).
$N$ represents some $k\in S/\mathcal W$, and by the
construction of $g$, $g(k)=t$. Thus, $g$ is onto.

If $f(s)=f(s')$, then the closures of $s$ and
$s'$ in the topology of $\mathcal W$ coincide. Since
the completion is Hausdorff, we must have
$\eta_{\mathcal W}(s)=\eta_{\mathcal W}(s')$;
thus, $g$ is one-one on the image of $\eta_{\mathcal W}$.
Suppose that $N$ and
$N'$ are Cauchy nets in
$\pair S{\mathcal W}$ which are not equivalent. Then
there is a $W\in\mathcal W$ such that $\pair
{N(d)}{N'(d')}\notin W$ for arbitrarily large values of
$d$ and $d'$. But $f^{-1}(V)\subseteq W$  for some
$V\in\mathcal V$, by the definition of $\mathcal
W=f^{-1}(\mathcal V)$. Thus, we have
$\pair{f(N(d))}{f(N'(d'))}\notin V$ for arbitrarily large
values of $d$ and $d'$. It follows that $f\circ N$ and
$f\circ N'$ are not equivalent. Thus, by theorem
\ref{T:Equivalence}(3), they have different limits in
$\pair T{\mathcal V}$. It follows that $g$ is one-one.

It remains to show that $g^{-1}$ is uniformly continuous.
We have $\eta^{-1}_{\mathcal W}(\mathcal W/\mathcal W)=\mathcal W$, by
(2), and since $f=g\eta_{\mathcal W}$ we have
\[
f^{-1}(g^{-1})^{-1}(\mathcal W/\mathcal W)=\eta^{-1}_{\mathcal W}(\mathcal W/\mathcal W)=\mathcal W=f^{-1}(\mathcal V);
\]
it follows by corollary \ref{T:Dense2} that
$(g^{-1})^{-1}(\mathcal W/\mathcal W)=\mathcal V$, whence
$g^{-1}$ is uniformly continuous.

(8): Let $f=\eta_{\mathcal V/\mathcal U}\eta_{\mathcal U}:
S\to(S/\mathcal U)/(\mathcal V/\mathcal U)$. Since $\eta_{\mathcal U}(S)$
is dense in $S/\mathcal U$, and $\eta_{\mathcal V/\mathcal U}(S/\mathcal U)$ is dense in $(S/\mathcal U)/(\mathcal V/\mathcal U)$, it follows
that $f(S)$ is dense in $C\pair{S/\mathcal U}{\mathcal V/\mathcal U}$.
Then by (7), there is an isomorphism $g_{\mathcal
U,\mathcal V}:C\pair S{\mathcal W}\to C\pair{S/\mathcal
U}{\mathcal V/\mathcal U}$, where $\mathcal
W=f^{-1}((\mathcal V/\mathcal U)/(\mathcal V/\mathcal
U))$, such that
$f=g_{\mathcal U,\mathcal V}\eta_{\mathcal W}$. However,
we have
$\mathcal W=f^{-1}((\mathcal V/\mathcal U)/(\mathcal
V/\mathcal U))=\eta^{-1}_{\mathcal U}(\mathcal V/\mathcal
U)=\mathcal V$.

(9): We have
\begin{align*}
\eta_{\mathcal V}^{-1}g_{\mathcal U,\mathcal
V}^{-1} ((\mathcal W/\mathcal U)/(\mathcal V/\mathcal U))
&=\eta^{-1}_{\mathcal U}\eta^{-1}_{\mathcal V
/\mathcal U}((\mathcal W/\mathcal U)/(\mathcal V
/\mathcal
U)) \\
 &=\eta^{-1}_{\mathcal U}(\mathcal W/\mathcal U) \\
&=\mathcal W,
\end{align*}
and also $\eta_{\mathcal V}^{-1}(\mathcal W/\mathcal V)=\mathcal W$ by
(2). Then since the image of $\eta_{\mathcal V}$ is dense, we
have the desired result by corollary~\ref{T:Dense2}.
 \end{proof}

\subsection{Hausdorff completion and uniform algebras}
 A uniform algebra $\pair A{\mathcal U}$ in
$\mathbf V$ is complete and Hausdorff iff it is an object of
$\mathbf V[\CHUnif]$.

\begin{theorem}
\label{T:AlgCompletion}
Let $\pair A{\mathcal U}\in\mathbf
V[\Unif]$. We have
\begin{enumerate}

\item[(1)] $C\pair A{\mathcal U}$ has a unique
structure of algebra in the variety $\mathbf V$ such that
$\eta_{\mathcal U}:\pair A{\mathcal U}\to C\pair
A{\mathcal U}$ is a homomorphism and $\pair A{\mathcal
U}$ is a uniform algebra;

\item[(2)] if $\mathcal V\in\OpUnif A$ is such
that $\mathcal U\leq\mathcal V$, then $\mathcal V/\mathcal U\in\OpUnif A/\mathcal U$;

\item[(3)] if $\mathcal V\in\OpUnif A$ is such that
$\mathcal U\leq\mathcal V$, then the isomorphism of
uniform spaces
$g_{\mathcal U,\mathcal V}$ of
theorem~\ref{T:Completion}(8) is an algebra isomorphism;

 \item[(4)] if $\alpha$ is a congruence
contained in $\mathcal U$, then $({\mathcal
U}\{\,\alpha\,\})/\mathcal U$ is a congruence on
$A/\mathcal U$;

\item[(5)] if $\mathcal V$ is a
congruential uniformity such that $\mathcal U\leq\mathcal V$, then
$\mathcal V/\mathcal U$ is congruential;

\item[(6)] $\pair{C\pair A{\mathcal U}}{\eta_{\mathcal
U}}$ is a universal arrow from $\pair A{\mathcal U}$ to
the forgetful functor from the category $\mathbf
V[\CHUnif]$ to
$\mathbf V[\Unif]$; and

\item[(7)] if $A$ is an algebra, $\pair B{\mathcal V}$
is a complete, Hausdorff uniform algebra, and $f:A\to B$
is a homomorphism with image dense in $B$, then $\pair
B{\mathcal V}$ is isomorphic to $C\pair A{f^{-1}(\mathcal
V)}$.

\end{enumerate}
\end{theorem}

\begin{proof} (1):
Let $t$ be an $n$-ary term, and let $\vec k$ be an
$n$-tuple of elements of $A/\mathcal U$. Let $N_i:\mathbf
D_i\to A$ be representatives of $k_i$ for each $i$. We
define $t(\vec N):\Pi_i\mathbf D_i\to S$ to be the
Cauchy net given by
\[
t(\vec N)(\vec d\,)=t(N_1(d_1),\ldots,N_n(d_n)),
\]
and $t(\vec k\,)$ to be the equivalence class of $t(\vec
N)$. Since $t:A^n\to A$ is uniformly continuous, this is
independent of the chosen representatives $N_i$. Thus we
have mapped $t$ to an $n$-ary operation on $A/\mathcal U$.

We claim that these mappings, for $n\geq 0$,
constitute a clone homomorphism. If $t=\pi_{in}$, then
$t(\vec N)(\vec d\,)=N_i(d_i)$; it is straightforward to
prove that
$t(\vec N)\sim N_i$. Thus, $\pi_{in}(\vec k)=k_i$.  On the
other hand, suppose
$t=v'\vec v$ where $v'$ is an $n'$-ary term and $\vec v$
is an
$n'$-tuple of $n$-ary terms. We have $v'(\vec v(\vec
N)):(\Pi_i\mathbf D_i)^{n'}\to S$, and $t(\vec
N)=v'(\vec v(\vec N))\circ\Delta$, where
$\Delta:\Pi_i\mathbf D_i\to(\Pi_i\mathbf D_i)^{n'}$ is
the diagonal map; it follows easily from this that
$t(\vec N)\sim v'(\vec v(\vec N))$, or $(v'\vec v)(\vec
k)=v'(\vec v(\vec k))$. Thus, we have a clone homomorphism
and
$A/\mathcal U$ is an algebra in
$\mathbf V$.

 It
is clear that $\eta_{\mathcal U}$ is a homomorphism. The
uniqueness of the algebra structure follows from the fact
that the image of $\eta_{\mathcal U}$ is dense in $C\pair
A{\mathcal U}$.

(2): Let $t$ be an $n$-ary term for $n>0$,
and let $V\in\mathcal V$. Then if $k_1\mathrel{R(\mathcal
U,V^t)}k_1'$, there exist $N_1\in k_1$, $N_1'\in k_1'$
such that
$N_1(d)\mathrel{V^t}N_1'(d')$ for large enough $d$ and
$d'$. Given any $k_2$, $\ldots$, $k_n\in A/\mathcal U$, with
representatives $N_2$, $\ldots$, $N_n$ respectively, we
have $t(N_1,N_2,\ldots,N_n)(\vec d)\mathrel
Vt(N_1',N_2,\ldots,N_n)(\vec d)$ for large enough $\vec
d$. This implies
\[
t(k_1,k_2,\ldots,k_n)\mathrel{R(\mathcal
U,V)}t(k_1',k_2,\ldots,k_n),
\]
 showing that $t^{A/\mathcal U}$ is uniformly  continuous
with respect to $\mathcal V/\mathcal U$ in the first
argument; by theorem \ref{T:Compatible}, $\mathcal
V/\mathcal U$ is compatible.

(3): If $\mathcal U\leq\mathcal V$, then by theorem
\ref{T:Completion}(8), the diagram of uniform spaces
\[
\begin{CD}
 \pair A{\mathcal U} @>\eta_{\mathcal U}>>
C\pair A{\mathcal U}
\\
@V\eta_{\mathcal V}VV @VV\eta_{\mathcal V/\mathcal U}V \\
C\pair A{\mathcal V} @>>g_{\mathcal U,\mathcal
V}>
C\pair{A/\mathcal U}{\mathcal V/\mathcal U} 
\end{CD}
\]
commutes. Since $\eta_{\mathcal U}$,
$\eta_{\mathcal V}$, and $\eta_{\mathcal V/\mathcal U}$
are homomorphisms by (1) and (2),
$g_{\mathcal U,\mathcal V}$ is a homomorphism on the
dense subalgebra
$\eta_{\mathcal V}(A)$ of $A/\mathcal V$. Since all of
the functions in the diagram are uniformly continuous, it
follows that
$g_{\mathcal U,\mathcal V}$ is an algebra homomorphism.

(4) is clear.

(5): Follows from (4) and theorem \ref{T:Base3}.

(6): It remains to prove only that if $f:\pair
A{\mathcal U}\to \pair B{\mathcal V}$ is a uniformly
continuous homomorphism and
$\pair B{\mathcal V}\in\mathbf V[\CHUnif]$, then the uniformly
continuous function $g:C\pair A{\mathcal U}\to\pair B{\mathcal V}$, given by the universal
property of the completion, is a homomorphism.
We have $f=g\eta_{\mathcal U}$, so the restriction
of $g$ to the dense subalgebra $\eta_{\mathcal U}(A)$ of $A/\mathcal U$ is a homomorphism. Since
$f$ and $g$ are uniformly continuous, it
follows that $g$ is a homomorphism.

(7): The universal property of $C\pair A{f^{-1}(\mathcal
V)}$ gives a uniformly continuous homomorphism $g:C\pair
A{f^{-1}(\mathcal V)}\to\pair B{\mathcal V}$, which is an
isomorphism of uniform spaces by theorem
\ref{T:Completion}(7). $g$ is a homomorphism because its
restriction to the dense subspace $\eta_{f^{-1}(\mathcal
V)}(A)$ is a homomorphism.
\end{proof}

Note that by (5), Hausdorff completion takes objects of
$\CU[\mathbf V]$ into $\CU[\mathbf V]$. We denote the category
of complete, Hausdorff, congruentially uniform algebras in
$\mathbf V$ by $\CHCU[\mathbf V]$.

\section{Hausdorff Completion of Congruentially
Uniform Algebras}\label{S:Completion}

In this section, we discuss a different formula for the
Hausdorff completion, which applies to the case of $\pair
A{\mathcal U}$ for $\mathcal U$ a congruential uniformity. This
definition is commonly used in commutative algebra; for
example, see \cite{A-M}.

Let $A$ be an algebra and $F\in\Fil\Con A$. We will
express the completion $A/({\mathcal U}F)$ as a
limit (also called an \emph{inverse limit}). We view $F$
as an inversely-directed set under inclusion. That is, we
will consider $F$ to be a category with objects
congruences
$\alpha\in F$ and with a single arrow from $\alpha$ to
$\beta\in F$ whenever $\alpha\leq\beta$. Let the functor
$D:F\to\mathbf V$ be defined by $D(\alpha)=A/\alpha$ and
for $f:\alpha\to\beta$, by letting $D(f):A/\alpha\to
A/\beta$ be the unique homomorphism such that
$D(f)\nat\alpha=\nat\beta$. Let $A_F=\varprojlim D$, with
limiting cone components $\xi_\alpha:A_F\to A/\alpha$,
and let $\eta_F:A\to A_F$ be the unique homomorphism
(given by the universal property of the limit) such that
$\xi_\alpha\eta_F=\nat\alpha$ for each $\alpha$.

\begin{theorem} There is an isomorphism
$\phi:A_F\cong A/({\mathcal U}F)$ such that
$\phi\eta_F=\eta_{\mathcal U}$.
\end{theorem}

\begin{proof} Each element $c\in A_F$ can be seen as a
choice of elements $c_\alpha\in A/\alpha$ for each
$\alpha$, such that whenever $\alpha\leq\beta$,
$c_\beta=c_\alpha/(\beta/\alpha)$. Given $c$, let $\bar
c(\alpha)\in A$ be some element such that $\bar
c(\alpha)/\alpha=c_\alpha$, for each $\alpha\in F$. Then
$\bar c$ is a Cauchy net, with respect to
${\mathcal U}F$, when $F$ is viewed as a
directed set given by the opposite of the
inversely-directed set described above. It thus represents
an element
$k\in A/({\mathcal U}F)$, and it is clear that this
element is independent of the choices made in the
definition of
$\bar c$. We define
$\phi(c)=k$.

It is straightforward to prove that $\phi$ is a
homomorphism.

To prove $\phi$ is onto, let $N\in k\in
A/({\mathcal U}F)$, where $N:\mathbf D\to A$. For
each $\alpha\in F$, there is an $\alpha$-class $c_\alpha$
such that $N(d)\in c_\alpha$ for large enough $d$. The
$c_\alpha$ determine an element $c$ of $A_F$ such that
$\phi(c)=k$.
\end{proof}

\section{Limits and Colimits in Categories of Uniform
Algebras}\label{S:Limits}

\subsection{Categories of uniform algebras}
For every
variety of algebras
$\mathbf V$, we have defined a number of categories,
which we arrange in the following diagram along with the
category $\mathbf V$:
\[
\begin{CD} \CHCU[\mathbf V]
@>>> \mathbf V[\CHUnif] \\
@VVV @VVV \\
\HCU[\mathbf V]
@>>> \mathbf V[\HUnif] \\
@VVV @VVV \\
\CU[\mathbf V]
@>>> \mathbf V[\Unif]
 \\
@. @VVV \\
@. \mathbf V
\end{CD}
\]

All of the arrows represent forgetful functors
which forget the uniformity or some property of it,
and all have left adjoints,
among them $H:\mathbf V[\Unif]\to\mathbf V[\HUnif]$, the
functor of Hausdorffization, and a functor $\bar
C:\mathbf V[\HUnif]\to\mathbf V[\CHUnif]$ such that
$\bar CH$ is naturally isomorphic to $C$, the functor of
Hausdorff completion. The forgetful functor from
$\mathbf V[\Unif]$ to $\mathbf V$ also has a right
adjoint.

\subsection{Limits and colimits in $\mathbf V[\Unif]$}
As the underlying algebra functor to $\mathbf V$ is both a
right adjoint functor and a left adjoint functor (see
theorem \ref{T:Adjoint}), it preserves all categorical
limits (sometimes called
\emph{inverse limits}), including products, and also all
colimits (sometimes called
\emph{direct limits}). Thus, to find a limit or colimit
in
$\mathbf V[\Unif]$, assuming it exists, we need only find
the corresponding limit or colimit in $\mathbf V$ and
provide it with the correct uniformity.

Limits of all small diagrams exist in
$\mathbf V[\Unif]$; the rule for finding a limit is to find
the limit of the underlying diagram in the category of
sets, and provide that limit with the weakest uniformity
such that all of the components of the cone to the diagram
(i.e., in the case of products, the projections) are
uniformly continuous. That weakest uniformity is the meet
of the inverse images $\xi_d^{-1}(\mathcal U_d)$, where $\xi_d$
is a component of the limiting cone, and $\mathcal U_d$ is the
uniformity for the corresponding object in the diagram.

Colimits of all small diagrams in $\mathbf V[\Unif]$
exist as well. The rule for finding the colimit of a
small diagram
$D:\mathbf D\to\mathbf V[\Unif]$ is to form the colimit $S$ of
the underlying diagram in $\mathbf V$, and then provide $S$
with the strongest uniformity such that all the
homomorphisms in the colimiting cone are uniformly
continuous. That uniformity is the join $\mathcal
U=\bigvee_{d\in\mathbf D}\mathcal U_d$, where $\mathcal
U_d$ is the direct image uniformity $\xi_d(\mathcal
U_d)$, with $\mathcal U_d$ being the uniformity of
$Dd$.

\subsection{Properties of $H$ and $C$}
Neither $H$ nor $C$ preserves pullbacks. However we
have
\begin{theorem} Both $H$ and $C$ preserve all products.
\end{theorem}

\subsection{Limits and colimits in the other categories}
In the other categories, the limit of a small diagram can
be found in $\mathbf V[\Unif]$ and is the limit in the
category in question because, as can be proved in each
case, it is a uniform algebra in that category. The
colimit of any small diagram can be found in $\mathbf
V[\Unif]$ and made into a uniform algebra in the category
in question by applying the left adjoint functor to the
appropriate forgetful functor; this results in the
colimit of the diagram in the category in question.

\section{Factorization Systems}\label{S:FactSys}

In this section, we review the standard
category-theoretic notion of a factorization system in a
category.

If $\mathbf C$ is a category, then a \emph{factorization
system} in $\mathbf C$ is a pair $\pair{\mathbf E}{\mathbf M}$
of subcategories of $\mathbf C$, having the following four
properties:
\begin{enumerate}
\item[(F1)] If $e$ is an arrow of $\mathbf E$ and $h$ is
an isomorphism such that the composite $he$ exists, then
$he\in\mathbf E$, while if $m$ is an arrow of $\mathbf M$ and
$h$ is an isomorphism such that the composite $mh$
exists, then $mh\in\mathbf M$.

\item[(F2)] Every arrow $f$
of $\mathbf C$ factors as $f=me$, with $m\in\mathbf M$ and
$e\in\mathbf E$.

\item[(F3)] Every commutative square
of the form
\[
\CD x @>>>
z \\
@VeVV @VVmV \\
y @>>> w\endCD
\]
with $e\in\mathbf E$ and $m\in\mathbf M$, admits a unique
diagonal arrow $\delta$ making the diagram
\[
\rlap{
\begin{picture}(20,20)
\put (10,-8) {\vector(4,3){30}}
\end{picture}
}
\rlap{\lower.06in\hbox{\qquad\ \ $\scriptstyle\delta$}}
\CD x @>>>
z \\
@VeVV @VVmV \\
y
@>>> w\endCD
\]
commutative.

\item[(F4)] Every isomorphism belongs both to $\mathbf E$
and to $\mathbf M$.
\end{enumerate}

These conditions are somewhat redundant; for example,
(F4) follows from the other three conditions. In any
case, if the conditions are satisfied, $\mathbf E$ and
$\mathbf M$ satisfy many additional properties as described
in, e.g., \cite{A-H-S}. Examples of this situation: (1)
the category of sets $\Set$, the subcategory $\mathbf E$ of
sets and onto functions, and the subcategory $\mathbf M$ of
sets and one-one functions; (2) more generally, onto
homomorphisms and one-one homomorphisms in the category of
all algebras in a variety of algebras; (3) the category of
small categories and functors, the subcategory $\mathbf E$
of \emph{cofaithful} functors, i.e., functors which are
one-one and onto on objects and onto each homset of the
codomain category, and the subcategory $\mathbf M$ of
faithful functors.

Two additional conditions that a factorization structure
$\pair{\mathbf E}{\mathbf M}$ can satisfy are that $\mathbf E$
consist entirely of epi arrows of $\mathbf C$, and that
$\mathbf M$ consist entirely of monic arrows. The
first two examples satisfy both conditions, and
cofaithful functors are epi as arrows of the category
of small categories and functors, but faithful
functors are not necessarily monic. 

\subsection{Subobject and quotient lattices}
If we have a class of epimorphisms $\mathbf E$ in a
category
$\mathbf C$, and $A\in\mathbf C$,
then for any two arrows $e:A\to B$, $e':A\to B'$ of $\mathbf
E$, there exists at most one arrow $f:B\to B'$ such that
$e'=fe$. If we then write $e\leq e'$ and identify $e$ and
$e'$ whenever both $e\leq e'$ and $e'\leq e$, we obtain a
partially-ordered set of equivalence classes associated
with $A$, which in many cases is a small set. (If this
holds for all $A\in\mathbf C$, then in a slight
extension of the usual terminology, we say $\mathbf E$ is
\emph{co-well-powered}.) We call this the $\mathbf
E$-quotient poset of $A$. We will see that under mild
conditions, when $\mathbf E$ is part of a
factorization system for $\mathbf C$, this
partially-ordered set is a lattice and we call it the
\emph{$\mathbf E$-quotient lattice} of $A$. For example,
starting from a set $A$ and the class $\mathbf E$ in
$\Set$ given by the onto functions, we obtain the lattice
of equivalence relations on $A$ as the $\mathbf
E$-quotient lattice of $A$; more generally, given an
algebra $A$ in a variety of algebras and letting $\mathbf
E$ be the subcategory of onto homomorphisms, we obtain
the congruence lattice, $\Con A$.

Dual considerations apply when we have a class of monic
arrows $\mathbf M$, and lead to a poset (small if $\mathbf M$
is \emph{well-powered}) of equivalence classes of arrows
of
$\mathbf M$ with codomain $A\in\mathbf C$, called the
$\mathbf M$-subobject poset of $A$.
Again, we call it the \emph{$\mathbf M$-subobject lattice}
if it is a lattice. In the same case of $\mathbf C=\Set$
(or, a variety of algebras) and $A$ a set (respectively,
an algebra), we obtain the lattice of subsets of the set
$A$ (respectively, the lattice of subalgebras).

Assume $\mathbf E$ and $\mathbf M$ 
form a factorization system $\pair{\mathbf E}{\mathbf M}$
in $\mathbf C$. In this case, it follows from the axioms
for a factorization system that if
$f$ is an arrow of
$\mathbf C$ and $f=me$ is a factorization, then the
equivalence classes of $m$ (in the $\mathbf M$-subobject
poset) and $e$ (in the $\mathbf E$-quotient poset) do not
depend on the particular factorization chosen, but only
on $f$.

\subsection{Factorization systems and lattice
operations}
Let $\mathbf C$ be a category and $\pair{\mathbf
E}{\mathbf M}$ a factorization system, such that
$\mathbf E$ consists of epi and $\mathbf M$ of monic
arrows.

\begin{proposition} We have
\begin{enumerate}

\item[(1)] If $\mathbf C$ has $n$-fold pushouts of
arrows in $\mathbf E$ for some cardinal $n$, then the
$\mathbf E$-quotient poset of each object of $\mathbf C$
has
$n$-fold joins; and

\item[(2)] if $\mathbf C$ has $n$-fold products for some
cardinal $n$, then the $\mathbf E$-quotient poset of each
object of $\mathbf C$ has $n$-fold meets.

\end{enumerate}
\end{proposition}

A dual statement gives sufficient
conditions for the $\mathbf M$-subobject posets to have
meets and joins.

\section{Factorization in $\mathbf
V[\CHUnif]$}\label{S:Factorization}

In this section, we will define an important factorization
system in the
category $\mathbf V[\CHUnif]$ of complete, Hausdorff uniform
algebras in $\mathbf V$.

Let $\mathbf M$ be the subcategory of uniformly continuous
homomorphisms $f:\pair A{\mathcal U}\to\pair B{\mathcal V}$ such
that $f$ is one-one and $\mathcal U=f^{-1}(\mathcal V)$.
Let $\mathbf E$ be the subcategory of uniformly continuous
homomorphisms $f:\pair A{\mathcal U}\to\pair B{\mathcal V}$ such
that the image of $f$ is dense in $B$.

\begin{theorem} $\pair{\mathbf E}{\mathbf M}$ is a
factorization system in $\mathbf V[\CHUnif]$.
\end{theorem}

\begin{proof}
It suffices to verify (F1) through (F3).

(F1): Follows from the fact that an isomorphism in $\mathbf
V[\CHUnif]$ is an algebra isomorphism $h:\pair
A{\mathcal U}\to\pair B{\mathcal V}$ such that  $\mathcal
U=h^{-1}(\mathcal V)$ and $\mathcal
V=(h^{-1})^{-1}(\mathcal U)$.

(F2): Given $f:\pair A{\mathcal U}\to \pair B{\mathcal V}$, let
$C$ be the closure of the image of $f$ in $\pair
B{\mathcal V}$, and $\mathcal W$ the inverse image
uniformity
$m^{-1}(\mathcal V)$ for $m:C\to B$ the inclusion.
$C$, being the closure of a subalgebra of $B$, is also
a subalgebra, and $\pair C{\mathcal W}\in\mathbf
V[\CHUnif]$. Then
$m\in\mathbf M$, and
$f$, considered as an arrow from $\pair A{\mathcal U}$ to
$\pair C{\mathcal W}$, is an arrow of $\mathbf E$.

(F3): Given a commutative square of the form
\[
\CD\pair A{\mathcal U} @>f>> \pair B{\mathcal V} \\
@VeVV @VVm,V \\
\pair C{\mathcal W} @>>g> \pair D{\mathcal Y} \endCD
\]
with $e\in\mathbf E$ and $m\in\mathbf M$, for any $c\in
C$, let
$N$ be a net in
$A$ such that
$e\circ N\rightarrow c$. We have $g\circ e\circ
N\rightarrow g(c)$.  Therefore, $m\circ f\circ
N\rightarrow g(c)$. However, $\mathcal V=m^{-1}(\mathcal
Y)$; it follows by theorem \ref{T:Cauchy}(3) that
$f\circ N$ is a Cauchy sequence with respect to $\mathcal
V$. Since
$\pair B{\mathcal V}$ is complete, it follows that
$g(c)\in B$. We define $\delta(c)=g(c)$. It is clear that
$\delta e=f$ and $m\delta=g$.

To show $\delta$ is uniformly continuous, let
$V\in\mathcal V$. Then $V=m^{-1}(Y)$ for some
$Y\in\mathcal Y$. Since
$m\delta=g$,
$\delta(g^{-1}(Y))=m^{-1}(Y)=V$. \end{proof}

Note that $\mathbf M$ is well-powered and $\mathbf E$ is
co-well-powered. Also, $\mathbf E$ consists of epi and
$\mathbf M$ of monic arrows.

\begin{theorem} The $\mathbf E$-quotient lattice of
a complete, Hausdorff uniform algebra $\pair A{\mathcal U}$
is isomorphic to $\I_{\OpUnif A}[\mathcal U,\top]$.
\end{theorem}

\begin{proof}
Given an arrow $f:\pair A{\mathcal U}\to\pair B{\mathcal
V}$, representing an element $[f]$ of the $\mathbf
E$-quotient lattice, let
$\phi([f])=f^{-1}(\mathcal V)$. $\phi$ is clearly
well-defined.
$\phi$ is onto $I_{\OpUnif A}[\mathcal U,\top]$ because
if $\mathcal W\in\OpUnif A$, with
$\mathcal U\leq\mathcal W$, then $\phi([\eta_{\mathcal
W}])=\mathcal W$.
$\phi$ is one-one because if
$\phi([f])=\mathcal W$ and $\phi([f'])=\mathcal W$, where
$f:\pair A{\mathcal U}\to\pair B{\mathcal V}$ and
$f':\pair A{\mathcal U}\to\pair{B'}{\mathcal V'}$, then
both $\pair B{\mathcal V}$ and $\pair{B'}{\mathcal V'}$
are isomorphic to $C\pair A{\mathcal W}$ by theorem
\ref{T:AlgCompletion}(7), by isomorphisms
$g:C\pair A{\mathcal W}\to\pair B{\mathcal V}$ and
$g':C\pair A{\mathcal W}\to\pair{B'}{\mathcal V'}$
satisfyin $f=g\eta_{\mathcal W}$ and
$f'=g'\eta_{\mathcal W}$, implying that
$f$ and
$f'$ represent the same element $[\eta_{\mathcal W}]$ of
the
$\mathbf E$-quotient lattice of $A$.
\end{proof}

\subsection*{Acknowledgement} We would like to thank
Keith Kearns and Michael Kinyon for their comments and
suggestions on this subject.

\end{document}